\documentclass{amsart}

\usepackage{amsmath, amssymb, mathrsfs, dsfont}

\numberwithin{equation}{section}

\begin{document}

\title[$L$-functions in positive characteristic]{On automorphic $L$-functions in \\ positive characteristic}


\author[L. A. Lomel\'i]{Luis Alberto Lomel\'i}


\date{May 2012}

\begin{abstract}
We give a proof of the existence of Asai, exterior square, and symmetric square local $L$-functions, $\gamma$-factors and root numbers in characteristic $p$, including the case of $p = 2$. Our study is made possible by developing the Langlands-Shahidi method over a global function field in the case of a Siegel Levi subgroup of a split classical group or a quasi-split unitary group. The resulting automorphic $L$-functions are shown to satisfy a rationality property and a functional equation. A uniqueness result of the author and G. Henniart allows us to show that the definitions provided in this article are in accordance with the local Langlands conjecture for ${\rm GL}_n$. Furthermore, in order to be self contained, we include a treatise of $L$-functions arising from maximal Levi subgroups of general linear groups.
\end{abstract}

\maketitle

\section*{Introduction}

We develop the Langlands-Shahidi method over a global function field in the case of a Siegel Levi subgroup of a split classical group or a quasi-split unitary group. In particular, we complete the study of exterior and symmetric square $L$-functions of \cite{hl11,lomeli2009} and establish the existence of Asai $L$-functions, and related local factors, in positive characteristic; these are uniquely characterized in \cite{hlarxiv}. The cases treated in this article include an induction step, necessary to develop the $\mathcal{LS}$ method for the classical groups in positive characteristic \cite{lomeliLSRH}.

A thorough understanding of automorphic $L$-functions and related local factors for the classical groups also requires Godement-Jacquet $L$-functions \cite{gj,j} as well as Rankin-Selberg products for representations of ${\rm GL}_m$ and ${\rm GL}_n$ \cite{jpss}. Indeed, Godement-Jacquet factors appear intrinsically in our study of the symplectic group and Rankin-Selberg products, which are ubiquitous in the Langlands program, are needed in a crucial multiplicativity property for representations obtained via unitary parabolic induction. In order to be self contained, our treatise includes $L$-functions arising from maximal Levi subgroups of general linear groups.

Thanks to a characterization of local factors with G. Henniart \cite{hl11,hlarxiv}, we can be certain that our definitions are in accordance with the Langlands conjectures for ${\rm GL}_n$ in positive characteristic \cite{laff,lrs}. 

While one purpose of the article is to complete our understanding of the Siegel Levi case for the split classical groups, the cases involving quasi-split unitary groups are new in positive characteristic. Thus, we establish existence of Asai $\gamma$-factors and prove that they satisfy a global functional equation involving partial $L$-functions; uniqueness of Asai $\gamma$-factors is proved in \cite{hlarxiv}. Furthermore, Asai $L$-functions are important in our study of automorphic $L$-functions in the case of a general maximal Levi subgroup of a unitary group \cite{lomeliU}. The following theorem summarizes the global results of this article (see Theorem~7.7); its proof begins on \S~1 and ends on \S~7.

\medskip

\noindent {\bf Theorem.} \emph{Let $k$ be a global function field with finite field of constants $\mathbb{F}_q$. Let ${\bf M}$ be a Siegel Levi subgroup of a split classical group or a quasi-split unitary group. Let $\pi$ be a cuspidal automorphic representation of ${\bf M}(\mathbb{A}_k)$. Then, the $L$-functions $L(s,\pi,r_i)$ are ${\rm nice}$, i.e., for each $i$:}
\begin{itemize}
   \item[(i)] (Rationality) \emph{$L(s,\pi,r_i)$ has a meromorphic continuation to a rational functions on $q^{-s}$;}
   \item[(ii)] (Functional equation) $L(s,\pi,r_i) = \varepsilon(s,\pi,r_i) L(1-s,\tilde{\pi},r_i)$.
\end{itemize}

\medskip

\noindent{\bf Remark:} \emph{Over a number field, global $L$-functions are said to be $\rm nice$ if they have a meromorphic continuation, are bounded on vertical strips, and satisfy a global functional equation. In the function field case, the rationality of $L(s,\pi,r_i)$ implies boundedness on vertical strips away from poles.}

\medskip

The automorphic $L$-functions of the theorem are indexed by representations $r_i$ depending on the classical group, which we now describe. Let $F$ be a non-archimedean local field of characteristic $p$ and let $\psi$ be a non-trivial character of $F$. Consider ${\bf M} \simeq {\rm GL}_n$ as a Siegel Levi subgroup of ${\bf G} = {\rm SO}_{2n+1}$ or ${\rm SO}_{2n}$. Let $r = {\rm Sym}^2$ when ${\bf G} = {\rm SO}_{2n+1}$, and let $r = \wedge^2$ when ${\bf G} = {\rm SO}_{2n}$. Given a smooth irreducible representation $\pi$ of ${\rm GL}_n(F)$, let $\sigma$ be the $\ell$-adic $n$-dimensional Frobenius semisimple representation of the Weil-Deligne group corresponding to $\pi$ under the local Langlands correspondence \cite{lrs}. The Langlands-Shahidi local coefficient $C_\psi(s,\pi,w_0)$, defined via intertwining operators of induced representations and the uniqueness property of Whittaker models when $\pi$ is $\psi$-generic, allows us to give a definition of $\gamma(s,\pi,r \circ \rho_n,\psi)$. In \cite{hl11}, the following equality is established
\begin{equation*}
   \gamma(s,\pi,r \circ \rho_n,\psi) = \gamma(s,r \circ \sigma,\psi),
\end{equation*}
where we have a Deligne-Langlands $\gamma$-factor on the right hand side. It is the case of a quasi-split even unitary group ${\rm U}_{2n}$ that leads us to the existence of Asai $\gamma$-factors. These are shown to be compatible with the local Langlands conjecture via tensor induction in \cite{hlarxiv}, where we establish the equality
\begin{equation*}
   \gamma(s,\pi,r_{\mathcal{A}},\psi) = \gamma(s,{}^\otimes{\rm I}(\sigma),\psi).
\end{equation*}
In all of these cases, we obtain only one $L$-function indexed by $r={\rm Sym}^2$, $\wedge^2$ or $r_\mathcal{A}$ depending on the quasi-split classical group. Next, when $\pi$ is $\psi$-generic, the cases of ${\rm Sp}_{2n}$ and ${\rm U}_{2n+1}$ yield
\begin{equation*}
   C_\psi(s,\pi,w_0) = \gamma(s,\pi,r_1,\psi) \gamma(2s,\pi,r_2,\psi),
\end{equation*}
where $\gamma(s,\pi,r_1,\psi)$ is a Godement-Jacquet $\gamma$-factor and $r_2$ is either an exterior square or a twisted Asai representation (see 6.4~Theorem').

We take the opportunity to build the theory from the ground up, beginning with abelian local factors of Tate's thesis, which arise as Langlands-Shahidi local coefficients of semisimple rank one classical groups. The theory of principal series is then presented; essential, since an automorphic representation is class one at almost every place. Local coefficients arise in the global theory via Fourier coefficients of appropriately chosen Eisenstein series and satisfy a crude functional equation. It is via $\gamma$-factors, defined by means of the local coefficient, that we define $L$-functions and root numbers. We remark that, recording results for future reference along the way, we are able to provide a thorough and linearly ordered proof of Theorem~7.7 beginning in \S~1. Also, we take \cite{conrad} to be an appropriate setting for working with quasi-split classical groups over function fields. In particular, unitary groups are defined over a degree-$2$ finite \'etale algebra $K$ over a field $k$.

Let us now give a more detailed description of the contents of the article. In \S~1, the necessary properties of $\gamma$-factors are gathered in the case of ${\rm GL}_1$ over a non-archimedean local field $F$. Already in this setting, local $L$-functions and root numbers, or $\varepsilon$-factors, are defined via $\gamma$-factors. In addition, given a degree-2 finite \'etale algebra $E$ over $F$, Langlands $\lambda$-factors give a connection between $\gamma$-factors over $E$ and $F$. The case of a separable quadratic algebra $E \simeq F \times F$ arises when considering quasi-split unitary groups over global fields at split places. In \S~2, we extend the basic definitions to this general setting. Recall that the local coefficient is defined via intertwining operators of parabolically induced representations and the uniqueness property of Whittaker models. Throughout the article, ${\bf G}$ is either a general linear group or one of the classical groups ${\rm SO}_{2n+1}$, ${\rm Sp}_{2n}$, ${\rm SO}_{2n}$, ${\rm U}_{2n}$ or ${\rm U}_{2n+1}$. When $\bf G$ is a classical group, the theory is studied for a standard Siegel parabolic subgroup ${\bf P} = {\bf M}{\bf N}$; the Siegel Levi subgroup $\bf M$ is isomorphic to ${\rm GL}_n$, except for the case of ${\rm U}_{2n+1}$, where it is isomorphic to ${\rm GL}_n \times {\rm U}_1$.

The semisimple rank one local computations of \S~3 are key to the Langlands-Shahidi method, they give compatibility of the method for the classical groups with local class field theory. The cases ${\bf G} = {\rm SL}_2$ and ${\rm SU}_3$ are mentioned in characteristic zero in \cite{kesh}, where it is shown how to extend these results via restriction of scalars. We provide a straight forward and self contained proof when ${\bf G} = {\rm U}_3$, including the case of residual characteristic two in the setting of a degree-$2$ finite \'etale algebra $E$ over $F$.

In \S~4, we fix a system of Weyl group element representatives and make the appropriate normalization of Haar measures. This is done via Shahidi's approach to Langlands' lemma, which gives an algorithm to decompose Weyl group elements, and the corresponding unipotent radical, leading towards a multiplicativity property of the local coefficient. Notice that we can extend the algorithm to include unitary groups over a degree-2 finite \'etale algebra $E$ over the local field $F$. When the parabolic subgroup is a Borel subgroup, multiplicativity allows one to reduce the theory of principal series to semisimple rank one cases; it can be seen as the formula of Gindikin and Karpelevi\v c in this setting.

The link to the global theory is made via the theory of Eisenstein series, which is available over a global function field \cite{harder,mw,morris}. We continue the discussion on the crude functional equation of \cite{lomeli2009} in \S~5 in order to include quasi-split unitary groups (see Theorem~5.1). For appropriately chosen Eisenstein series corresponding to globally generic cuspidal automorphic representations, the global intertwining operator appearing in their functional equation decomposes into a product of local intertwining operators discussed in \S~2. Whittaker models allow one to compute Fourier coefficients, making the connection to the local coefficient possible. The formula of Casselman-Shalika for Whittaker models of unramified principal series and a formula for the intertwining operator in terms of $L$-functions, when it is evaluated at the identity, provide the necessary local results required in the proof of the crude functional equation.

A refinement of the crude functional equation leads to a definition of $\gamma$-factors, given in \S~6. When ${\bf G} = {\rm GL}_{m+n}$, ${\rm SO}_{2n+1}$, ${\rm SO}_{2n}$ and ${\rm U}_{2n}$, the crude functional equation satisfied by the local coefficient involves only one partial $L$-function of a cuspidal automorphic representation of ${\bf M}(\mathbb{A}_k)$. These cases give the main induction of the Langlands-Shahidi method for the classical groups. When ${\bf G}={\rm Sp}_{2n}$, $n>1$, or ${\rm U}_{2n+1}$, the crude functional equation gives rise to two individual functional equations in terms of $\gamma$-factors. In \S~6.5, we provide all of the properties necessary to uniquely characterize $\gamma$-factors (see \cite{hl11,hlarxiv}). In the appendix, we give the characterization for the case of a general linear group.

We define local $L$-functions and root numbers, or $\varepsilon$-factors, in \S~7. First, for tempered representations. Then, in general via Langlands classification and analytic continuation. The global functional equation takes its final form in terms of completed $L$-functions and global root numbers in \S~7.6. Let $\pi$ be a cuspidal automorphic representation of ${\bf M}(\mathbb{A}_k)$, then
\begin{equation*}
   L(s,\pi,r_i) = \varepsilon(s,\pi,r_i) L(1-s,\tilde{\pi},r_i),
\end{equation*}
where $\tilde{\pi}$ denotes the contragredient of $\pi$ and the representations $r_i$ are described in \S~5.2. This includes exterior square, symmetric square, Asai and Rankin-Selberg $L$-functions. Notice that the functional equation of Tate's thesis over a global function field is obtained by studying a Gr\"ossencharakter viewed as an automorphic representation of the maximal torus ${\bf M} \simeq {\rm GL}_1$ of ${\bf G} = {\rm SL}_2$.

The approach taken in this article does not make use of Lie algebras, other than to check that our definition of $L$-functions agrees with the definition using $L$-groups via the adjoint action of ${}^LM$ on ${}^L\mathfrak{n}$ and the Satake parametrization for unramified principal series. For a general smooth irreducible representation, the results of \cite{hl11,hlarxiv} and the appendix, show that our $\gamma$-factors and related local factors agree with those of the corresponding representation of the Weil-Deligne group via the local Langlands correspondence. 

I would like to thank F. Shahidi for his guidance and support, in addition to many interesting mathematical discussions held over the course of several years. I also wanted to take the opportunity to thank G. Henniart for many insightful remarks; joint work on the characterization of $\gamma$-factors helped solidify the induction step of the Langlands-Shahidi method for the classical groups. Mathematical conversations with L. Lafforgue and J.-K. Yu were very helpful while writing this paper; I am thankful to them for taking interest in this project. Work was begun at Purdue University and continued at the Institut des Hautes \'Etudes Scientifiques during the summer of 2011, where the main part of this article was written. I am very grateful to these institutions for their hospitality and support.

\section{Abelian $\gamma$-factors}

Let $\mathcal{O}_F$ be the ring of integers of a non-archimedean local field $F$. Let $\mathfrak{p}_F$ be the maximal ideal of $\mathcal{O}_F$ and let $q_F$ denote the cardinality of $\mathcal{O}_F/\mathfrak{p}_F$. For every continuous non-trivial character $\psi : F \rightarrow \mathbb{C}^\times$ of level $l$, fix a Haar measure $\mu_\psi$ of $F$ such that $\mu_\psi(\mathcal{O}_F) = q_F^{l/2}$. Let $C_c^\infty(F)$ be the Bruhat-Schwartz space of complex valued locally constant functions on $F$ with compact support. The Fourier transform $\mathscr{F}_\psi : C_c^\infty(F) \rightarrow C_c^\infty(F)$ is defined by
\begin{equation*}
   \mathscr{F}_\psi(f)(x) = \int_F f(y) \psi(xy) d\mu_\psi(y).
\end{equation*}
Given a function $f$ with domain $F$ and $a \in F^\times$, let $f^a$ denote the function on $F$ defined by the rule $f^a(x) = f(ax)$. Then
\begin{equation}
   \mathscr{F}_{\psi^a} \circ \mathscr{F}_{\psi^b}(f) = \left| ab^{-1} \right|_F^{\frac{1}{2}} f^{-ab^{-1}},
\end{equation}
for every $f \in C_c^\infty(F)$, $a,b \in F^\times$. Notice that the level of $\psi^a$ is $l - {\rm ord}_F(a)$ and $\mu_{\psi^a} = \left| a \right|_F^{1/2} \mu_\psi$. Also, every continuous non-trivial additive character of $F$ is of the form $\psi^a$, for some $a \in F^\times$.

Fix Haar measures $\mu_\psi^\times$ on $F^\times$ so that
\begin{equation*}
   d\mu_\psi^\times(x) = \dfrac{q_F}{q_F - 1} \dfrac{d\mu_\psi(x)}{\left| x \right|_F}.
\end{equation*}
Then $\mu_\psi^\times(\mathcal{O}_F^\times) = \mu_\psi(\mathcal{O}_F)$.

Given a continuous character $\chi : F^\times \rightarrow \mathbb{C}^\times$, $f \in C_c^\infty(F)$, and $s \in \mathbb{C}$, the integral
\begin{equation*}
   \zeta_\psi(s,\chi,f) = \int_{F^\times} f(x) \chi(x) \left| x \right|_F^s d\mu_\psi^\times(x),
\end{equation*}
defines a Laurent power series on $q_F^{-s}$. These $\zeta$-functions converge and define a rational function on $q_F^{-s}$ in some right half plane ($\Re(s) > 0$ if $\chi$ is unitary). Given a continuous multiplicative character $\chi$, a continuous non-trivial additive character $\psi$, and $a \in F^\times$, the abelian $\gamma$-factor $\gamma(s,\chi,\psi^a)$ is the unique function on the variable $s$ satisfying Tate's local functional equation
\begin{equation}
   \zeta_\psi(1-s,\chi^{-1},\mathscr{F}_{\psi^a}(f)) = \gamma(s,\chi,\psi^a) \zeta_\psi(s,\chi,f),
\end{equation}
for every $f \in C_c^\infty(F)$. This $\gamma$-factor is a rational function on $q_F^{-s}$. Equation~(1.2), applied twice in combination with~(1.1), gives
\begin{equation}
   \gamma(s,\chi,\psi) \gamma(1-s,\chi^{-1},\overline{\psi}) = 1.
\end{equation}
The relationship between $\mu_{\psi^a}$ and $\mu_\psi$, together with~(1.2), yields
\begin{equation}
   \gamma(s,\chi,\psi^a) = \chi(a) \left| a \right|_F^{s - \frac{1}{2}} \gamma(s,\chi,\psi).
\end{equation}
Root numbers, or $\varepsilon$-factors, are defined via the equation
\begin{equation}
   \varepsilon(s,\chi,\psi) = \gamma(s,\chi,\psi) \dfrac{L(s,\chi)}{L(1-s,\chi^{-1})},
\end{equation}
where $L(s,\chi) = 1$ if $\chi$ is ramified, and $L(s,\chi) = ( 1 - q_F^{-s} \chi(\varpi_F))^{-1}$ if $\chi$ is unramified. They satisfy
\begin{subequations}
\begin{align}
   \varepsilon(s,\chi,\psi) & \varepsilon(1-s,\chi^{-1},\overline{\psi}) = 1, \\
   \varepsilon(s,\chi,\psi^a) & = \chi(a) \left| a \right|_F^{s-\frac{1}{2}} \varepsilon(s,\chi,\psi).
\end{align}
\end{subequations}
Additionally, if $\chi$ is unramified and $l = 0$ then $\varepsilon(s,\chi,\psi) = 1$, i.e.,
\begin{equation} \label{1unramified}
   \gamma(s,\chi,\psi) = \dfrac{L(1-s,\chi^{-1})}{L(s,\chi)}.
\end{equation}

Extend complex valued locally constant functions $f \in C^\infty(F^\times)$ to locally constant functions in $C^\infty(F)$ by setting $f(0) = 0$. Then $C_c^\infty(F) = {\rm span}(\mathds{1}_{\mathcal{O}_F}, C_c^\infty(F^\times))$, where $\mathds{1}_A$ denotes the characteristic function of a set $A \subset F$. Taking $N$ with $\chi_{1+\mathfrak{p}_F^N} = 1$ and setting $f = \mathds{1}_{1+\mathfrak{p}_F^N}$ in (1.2) yields the formula
\begin{equation*}
   \gamma(s,\chi,\psi) = \int_{\mathfrak{p}_F^{l-N}} \chi^{-1}(x) \left| x \right|_F^{-s} \psi(x) d\mu_\psi(x).
\end{equation*}
Here $\mathds{1}_{\mathfrak{p}_F^{l-N}}(x)\chi^{-1}(x) \left| x \right|_F^{-s} \in C_c^\infty(F^\times) \subset C_c^\infty(F)$. Therefore, the formula is valid for all $s \in \mathbb{C}$, and it is also possible to write
\begin{equation*}
   \gamma(s,\chi,\psi) = \int_F \chi^{-1}(x) \left| x \right|_F^{-s} \psi(x) d\mu_\psi(x),
\end{equation*}
which converges as a principal value integral. Combining this with~(3) gives
\begin{align*}
   \gamma(s,\chi,\psi)^{-1} & = \int_F \chi(x) \left| x \right|_F^{s-1} \overline{\psi}(x) d\mu_\psi(x) \\
   					  & = (1-q_F^{-1}) \int_{F^\times} \chi(x) \left| x \right|_F^s \overline{\psi}(x) d\mu_\psi^\times(x).
\end{align*}

Given a separable quadratic extension $E/F$ of locally compact non-archimedean fields, let $\eta_{E/F} : F^\times \rightarrow \mathbb{C}^\times$ be the character given by the rule: $\eta_{E/F}(x) = 1$ if $x \in {\rm N}_{E/F}(E^\times)$; $\eta_{E/F}(x) = -1$ if $x \notin {\rm N}_{E/F}(E^\times)$. Define
\begin{equation*}
   \lambda(E/F,\psi) = \dfrac{\int_{F^\times} \eta_{E/F}(x) \psi(x) d\mu_\psi^\times(x)}{\left| \int_{F^\times} \eta_{E/F}(x) \psi(x) d\mu_\psi^\times(x) \right|}.
\end{equation*}
It satisfies
\begin{equation*}
   \lambda(E/F,\psi) \lambda(E/F,\overline{\psi}) = 1.
\end{equation*}
To avoid confusion, given a continuous character $\chi': E^\times \rightarrow \mathbb{C}^\times$ and a continuous non-trivial character $\psi': E \rightarrow \mathbb{C}^\times$, write $\gamma_E(s,\chi',\psi')$ for the corresponding $\gamma$-factor. Then
\begin{equation} \label{langlandslambda}
   \lambda(E/F,\psi) \gamma_E(s,\chi \circ {\rm N}_{E/F},\psi \circ {\rm Tr}_{E/F}) = \gamma(s,\chi \vert_{F^\times},\psi) \gamma(s,\eta_{E/F} \chi \vert_{F^\times},\psi).
\end{equation}
To simplify notation, it is convenient to write $\psi_E = \psi \circ {\rm Tr}_{E/F}$.

Now, let $E \simeq F \times F$. We can assume $E$ is the separable quadratic algebra $F \times F$ by fixing a basis. For every $x = (x_1, x_2) \in E$, let $\left| x \right|_E = \left| x_1 x_2 \right|_F$, $\bar{x} = (x_2, x_1)$, ${\rm N}_{E/F}(x) = x \bar{x}$, and ${\rm Tr}_{E/F}(x) = x + \bar{x}$. A smooth character $\chi: E^\times \rightarrow \mathbb{C}^\times$ is of the form $\chi = \chi_1 \otimes \chi_2$, with $\chi_1$ and $\chi_2$ continuous characters of $F^\times$. Complex valued locally constant functions $f \in C^{\infty}(E^\times)$ are extended to complex valued locally constant functions $C^\infty(E)$ by setting $f(x)=0$ for $x \in E - E^\times$. Embed $F$ diagonally in $E$, in order to interpret $\chi \vert_{F^\times}$ and $\psi_E = \psi \circ {\rm Tr}_{E/F}$ correctly. Let $\mu_{\psi_E} = \mu_\psi \times \mu_\psi$. Then, equation~\eqref{langlandslambda} also holds in this case with $\eta_{E/F} = 1$ and $\lambda(E/F,\psi) = 1$, i.e.,
\begin{equation*}
   \gamma_E(s,\chi \circ {\rm N}_{E/F},\psi \circ {\rm Tr}_{E/F}) = \gamma(s,\chi_1 \chi_2, \psi)^2.
\end{equation*}

\section{On the local coefficient and the classical groups}

Given an algebraic group $\bf G$ defined over a field $k$, or an algebra $A$ over $k$, $G$ will denote the group of rational points, i.e., $G = {\bf G}(k)$ or ${\bf G}(A)$. Let $F$ be a non-archimedean local field.

\subsection{The general linear group} Let ${\bf G}$ be the general linear group ${\rm GL}_{n_1+n_2}$. Let $E$ be either the field $F$ or a degree-$2$ finite \'etale algebra over $F$. Let ${\bf B} = {\bf TU}$ be the Borel subgroup of $\bf G$ with maximal torus $\bf T$ and unipotent radical $\bf U$, whose group of rational points $B$ over $E$ consists of upper triangular matrices. The non-trivial character $\psi: F \rightarrow \mathbb{C}^\times$ extends to a character $U$, also denoted by $\psi$, by setting $\psi(u) = \psi_E(u_{1,2} + \cdots + u_{n_1+n_2-1,n_1+n_2})$, for $u = (u_{i,j}) \in U$. Here, $\psi_E = \psi \circ {\rm Tr}_{E/F}$ if $E = F \times F$, and $\psi_E = \psi$ if $E = F$. Consider the maximal Levi subgroup ${\bf M} \simeq {\rm GL}_{n_1} \times {\rm GL}_{n_2}$, whose group of rational points over $E$ is
\begin{equation*}
   M = \left\{ \left( \begin{array}{cc} g_1 & 0 \\ 0 & g_2 \end{array} \right) \vert \, g_1 \in {\rm GL}_{n_1}(E), g_2 \in {\rm GL}_{n_2}(E) \right\}.
\end{equation*}
Let ${\bf P} = {\bf MN}$ denote the standard parabolic subgroup of $\bf G$ with unipotent radical $\bf N$. 

Given smooth irreducible representations $\pi_1$ of ${\rm GL}_{n_1}(E)$ and $\pi_2$ of ${\rm GL}_{n_2}(E)$, let $\pi = \pi_1 \otimes \tilde{\pi}_2$. Here, $\tilde{\pi}_2$ denotes the contragredient representation of $\pi_2$. Consider the unitarily induced representation
\begin{equation*}
   {\rm I}(s,\pi) = {\rm ind}_{P}^G(\left| {\rm det}(\cdot) \right|_E^{\frac{s}{2}} \pi_1 \otimes \left| {\rm det}(\cdot) \right|_E^{-\frac{s}{2}} \tilde{\pi}_2).
\end{equation*}
More precisely, let ${\rm V}(s,\pi)$ be the space of functions $f: G \rightarrow V$ satisfying the following properties: there exists a compact open subgroup $K$ of $G$ such that $f(gk) = f(g)$ for every $g \in G$, $k \in K$; and, $f(mng)= \left| {\rm det}(m_1) \right|_E^{\frac{s}{2}}\pi_1(m_1) \otimes \left| {\rm det}(m_2)\right|_E^{-\frac{s}{2}} \tilde{\pi}_2(m_2) \delta_{P}^{\frac{1}{2}}(m) f(g)$, for every $m = (m_1,m_2) \in M$, $n \in N$, $g \in G$, where $\delta_P$ denotes the modulus character of $P$.  Then, ${\rm I}(s,\pi)$ is the right regular representation of $G$ on the space of functions ${\rm V}(s,\pi)$.

Elements of the Weyl group $W = W({\bf T},{\bf G)}$ will be identified with a fixed representative in ${\bf N}(F)$, where $F$ is embedded diagonally if $E = F \times F$. This is done globally in such a way that the theory matches the semisimple rank one local theory of \S~3, and similarly for Haar measures.  Let $w_0 = w_l w_{l,{\bf M}}$ the Weyl group element which is the product of the longest Weyl group element $w_l$ and the longest Weyl group element $w_{l,{\bf M}}$ with respect to $\bf M$. Let ${\bf M}' = w_0{\bf M}w_0^{-1}$, and let ${\bf P}'$ be the corresponding standard parabolic subgroup of $\bf G$ with unipotent radical ${\bf N}'$. Given a representation $\pi$ of $M$, let $w_0 (\pi)$ denote the representation of $M'$ given by $\pi(w_0^{-1} m' w_0)$, $m' \in M'$. It is $\psi_{M'}$-generic, where $\psi_{M'}(u') = \psi_M(w_0^{-1} u' w_0)$, $u' \in M' \cap U$. 

The intertwining operator ${\rm A}(s,\pi,w_0): {\rm V}(s,\pi) \rightarrow {\rm V}(-s,w_0(\pi))$, is defined by
\begin{equation*}
   {\rm A}(s,\pi,w_0)f(g) = \int_{N'} f(w_0^{-1}n'g) dn',
\end{equation*}
where the integral converges in a right half plane. In fact, it is a rational operator in the sense of Waldspurger, \S~IV of \cite{w}.

Let $C_c^\infty(G,V)$ be the space of locally constant functions with values in the space $V$ of $\pi$. For every $s \in \mathbb{C}$, there is a surjective map $\mathscr{P}_s : C_c^\infty(G,V) \rightarrow {\rm V}(s,\pi)$, given by
\begin{equation*}
   \mathscr{P}_s \varphi(g) = \int_{P}  \left| {\rm det}(m_{p,1}) \right|_E^{-\frac{s}{2}} \left| {\rm det}(m_{p,2}) \right|_E^{\frac{s}{2}} \delta_P^{\frac{1}{2}}(m_p) \pi^{-1}(m_p) \varphi(pg) dp.
\end{equation*}
Here, every $p \in P$ can be written as $p = m_p n_p$ for some $m_p \in M$, $n_p \in N$; furthermore, $m_p$ is identified with $(m_{p,1},m_{p,2}) \in {\rm GL}_{n_1}(E) \times {\rm GL}_{n_2}(E)$ by means of the embedding of $M$ inside $G$. 

Let $\mathcal{K}$ be a special maximal compact open subgroup of $G$. For each $s \in \mathbb{C}$, consider the function $f_s^0 = \mathscr{P}_s (\mathds{1}_{Pw_0B \cap \mathcal{K}})$. Notice that $\psi$ restricted to $M \cap U$ defines a character, which we denote by $\psi_M$. Then, if $\pi$ is $\psi_M$-generic with Whittaker functional $\lambda_{\psi_M}$, the function
\begin{equation*}
   \lambda_\psi(s,\pi) f_s^0 = \int_{N'} \lambda_{\psi_M} (f_s^0(w_0^{-1}n))\overline{\psi}(n) dn
\end{equation*}
defines a non-zero polynomial function on $\left\{ q_F^{-s}, q_F^{\,s} \right\}$ (Theorem~1.4 of \cite{lomeli2009}). It is possible to use the previous integral to define a Whittaker functional for ${\rm I}(s,\pi)$. The local coefficient is defined using the uniqueness property of Whittaker models and is given by
\begin{equation*}
   C_{\psi}(s,\pi,w_0) = \dfrac{\lambda_\psi(s,\pi)f_s^0}{\lambda_\psi(-s,w_0(\pi)) {\rm A}(s,\pi,w_0)f_s^0}.
\end{equation*}
It is a non-zero rational function on $q_F^{-s}$ \cite{lomeli2009,sha81}.

\subsection{Siegel Levi subgroups of classical groups}
Let ${\bf G} = {\bf G}_n$ be either a split classical group ${\rm SO}_{2n+1}$, ${\rm Sp}_{2n}$, ${\rm SO}_{2n}$, or a quasi-split unitary group ${\rm U}_{2n}$, ${\rm U}_{2n+1}$. A definition of the split classical groups over a non-archimedean local field or a global function field $k$ can be found in \S~3 of \cite{lomeli2009}, including the case ${\rm char} \, k = 2$. A definition of the quasi-split unitary groups, suitable for computing with local coefficients, is provided below.

Let $K$ be a degree-$2$ finite \'etale algebra over a field $k$. Let $\theta$ denote the non-trivial involution of $K$ over $k$. Write $\bar{x} = \theta(x)$, for $x \in K$. Extend conjugation to elements $g=(g_{i,j}) \in {\rm GL}_{m}(K)$, by setting $\bar{g}=(\bar{g}_{i,j})$. Given a positive integer $n$, consider
\begin{align*}
   h_{2n}(x) & = \sum_{i=1}^n \bar{x}_i x_{2n+1-i} - \sum_{i=1}^n \bar{x}_{2n+1-i} x_i, \ x \in K^{2n}, \\
   h_{2n+1}(x) & = \sum_{i=1}^{2n} \bar{x}_i x_{2n+2-i} - \bar{x}_{n+1} x_{n+1}, \ x \in K^{2n+1}.
\end{align*}
The unitary group ${\bf G} = {\rm U}_m$, $m = 2n$ or $m = 2n+1$, is defined so that its group of rational points is given by
\begin{equation*}
   G = \left\{ g \in {\rm GL}_m(K) \vert h_m(gx) = h_m(x) \text{ for every } x \in K^m \right\}.
\end{equation*}

\noindent{\bf Remark.} Throughout the article, we assume that split classical groups are defined over $k$ and non-split quasi-split classical groups are obtained via a degree-$2$ extension $K/k$.

For every positive integer $n$, let $\Phi_n$ be the $n \times n$ matrix with $ij$-entries $(\delta_{i,n-j+1})$. Fix the Siegel Levi subgroup $\bf M$ of ${\bf G}$ so that its group of rational points $M$ is embedded in $G$ as follows:
\begin{align*}
   M & = \left\{ \left( \begin{array}{ccc} g &  &  \\  & 1 &  \\  &  & \Phi_n {}^tg^{-1} \Phi_n \end{array} \right) \vert \, g \in {\rm GL}_{n}(k) \right\}, \text{ if } {\bf G} = {\rm SO}_{2n+1}; \\
   M & = \left\{ \left( \begin{array}{cc} g &  \\  & \Phi_n {}^tg^{-1} \Phi_n \end{array} \right) \vert \, g \in {\rm GL}_n(k) \right\}, \text{ if } {\bf G} = {\rm Sp}_{2n} \text{ or } {\rm SO}_{2n}; \\
   M & = \left\{ \left( \begin{array}{ccc} g & & \\ & z & \\ & & \Phi_n {}^t \bar{g}^{-1} \Phi_n  \end{array} \right) \vert \, g \in {\rm GL}_{n}(K), \, z \in {\rm ker}({\rm N}_{K/k}) \right\} \text{ if } {\bf G} = {\rm U}_{2n+1};\\
   M & = \left\{ \left( \begin{array}{cc} g & \\ & \Phi_n {}^t\bar{g}^{-1} \Phi_n \end{array} \right) \vert \, g \in {\rm GL}_n(K) \right\} \text{ if } {\bf G} = {\rm U}_{2n}.
\end{align*}
Let $\bf T$ be the maximal torus whose group of rational points $T$ consists of diagonal matrices, and let ${\bf B} = {\bf TU}$ be the Borel group of ${\bf G}$ with unipotent radical $\bf U$, whose group of rational points $B$ consists of upper triangular matrices. Let ${\bf P} = {\bf M}{\bf N}$ be the standard parabolic subgroup of ${\bf G}$ with Levi $\bf M$ and unipotent radical $\bf N$.

\subsection{} If ${\bf G}$ is a split classical group, let $E$ be the non-archimedean local field $F$. If ${\bf G}$ is a quasi-split unitary group, let $E$ be a degree-$2$ finite \'etale algebra over $F$. Hence $M \simeq {\rm GL}_n(E)$, unless ${\bf G} = {\rm U}_{2n+1}$. In this latter case, $M \simeq {\rm GL}_n(E) \times E^1$, where $E^1 = {\rm ker}({\rm N}_{E/F})$. Given a smooth irreducible representation $\pi$ of $M$, consider the unitarily induced representation
\begin{equation} \label{eq:ind}
   {\rm I}(s,\pi) = \left\{ \begin{array}{ll} {\rm ind}_P^G(\left| {\rm det}(\cdot) \right|_E^{s} \pi) & \text{ if } {\bf G} = {\rm SO}_{2n}, {\rm Sp}_{2n} \text{ or } {\rm U}_{2n+1} \\
                                   			  {\rm ind}_P^G(\left| {\rm det}(\cdot) \right|_E^{\frac{s}{2}}\pi) & \text{ if } {\bf G} = {\rm SO}_{2n+1} \text{ or } {\rm U}_{2n}
                                   \end{array} . \right.
\end{equation}
In the case of ${\rm U}_{2n+1}$, notice that a smooth irreducible representation $\pi$ of $M$ is of the form $\pi = \pi' \otimes \nu$, where $\pi'$ is a smooth irreducible representation of ${\rm GL}_n(E)$ and $\nu$ is a smooth character of $E^1$. Also in this case, the character $\left| {\rm det}(\cdot) \right|_E$ of ${\rm GL}_n(E)$ is extended to a character of $M$, trivial on $E^1$.

The space of ${\rm I}(s,\pi)$ is denoted by ${\rm V}(s,\pi)$. The Weyl group element $w_0 = w_l w_{l,{\bf M}}$ will be identified with a fixed representative in $G$. Notice that $\bf P$ is self-associate. The intertwining operator ${\rm A}(s,\pi,w_0): {\rm V}(s,\pi) \rightarrow {\rm V}(-s,w_0(\pi))$ is then defined by
\begin{equation*}
   {\rm A}(s,\pi,w_0)f(g) = \int_N f(w_0^{-1}ng)dn,
\end{equation*}
There is a surjective map $\mathscr{P}_s : C_c^\infty(G,V) \rightarrow {\rm V}(s,\pi)$, $s \in \mathbb{C}$, defined by
\begin{equation*}
   \mathscr{P}_s \varphi(g) = \int_P \left| {\rm det}(m_p) \right|_E^{-is} \delta_P^{\frac{1}{2}}(m_p) {\pi}^{-1}(m_p) \varphi(pg) dp,
\end{equation*}
where $i=1$ or $1/2$ depending on~(\ref{eq:ind}), and each $p \in P$ has $p  = m_p n_p$, $m_p \in M$, $n_p \in N$. Let $\mathcal{K}$ be a maximal compact open subgroup of $G$ such that $G = P\mathcal{K}$, and consider the function $f_s^0 = \mathscr{P}_s(\mathds{1}_{Pw_0B \cap \mathcal{K}})$. 

The non-trivial character $\psi: F \rightarrow \mathbb{C}^\times$ extends to a character of $\psi_E$ of $E$. As before, it then gives a character of the unipotent upper triangular matrices in ${\rm GL}_n(E)$. This defines a non-degenerate character $\psi_M$ of $M \cap U$. Let ${\bf G}_1$ be a classical group of semisimple rank one, and assume ${\bf G} = {\bf G}_n$ is a classical group of the same type as ${\bf G}_1$ with $n>1$. Given a non-negative integer $m$, let $I_m$ denote the $m \times m$ identity matrix. Then, the group of rational points $G_1$ is identified with a subgroup of $G$:
\begin{equation} \label{G_1embed}
   G_1 \simeq \left\{ \left( \begin{array}{ccc} I_{n-1} & & \\  & g & \\  & & I_{n-1} \end{array} \right) \vert \ g \in G_1 \right\} \subset G.
\end{equation}
The character $\psi_M$ is then extended to a character $\psi$ of $U$ so that $\psi \vert_{G_1 \cap U}$ agrees with the semisimple rank one definitions of the next section. Then, if $\pi$ is a $\psi_M$-generic irreducible representation of $M$ with Wittaker functional $\lambda_{\psi_M}$, the function
\begin{equation*}
   \lambda_\psi(s,\pi) f_s^0 = \int_N \lambda_{\psi_M}(f_s^0(w_0^{-1}n)) \overline{\psi}(n) dn,
\end{equation*}
is an non-zero Laurent polynomial function on $q_F^{-s}$ (Theorem~1.4 of \cite{lomeli2009}). The local coefficient is then given by
\begin{equation*}
   C_\psi(s,\pi,w_0) = \dfrac{\lambda_\psi(s,\pi) f_s^0}{\lambda_\psi(-s,w_0(\pi)){\rm A}(s,\pi,w_0)f_s^0}.
\end{equation*}
It is a non-zero rational function on $q_F^{-s}$ \cite{lomeli2009,sha81}.

\section{The local coefficient for semisimple rank one classical groups}

\subsection{Definitions} First, let ${\bf G} = {\rm U}_2$ or ${\rm U}_3$, let $E$ be either the algebra $F \times F$ or let $E/F$ be a separable quadratic extension of non-archimedean local fields. Let $E \rightarrow E$, $x \mapsto \bar{x}$, denote conjugation on $E$. If $E = F \times F$, let $\beta = (0,1)$. Assume that $E/F$ is an extension of local fields, then: if ${\rm char}(F) =2$, there is a $\beta \in \mathcal{O}_E^\times$ with galois conjugate $\bar{\beta} \in \mathcal{O}_E^\times$ such that ${\rm Tr}_{E/F}(\beta) \in \mathcal{O}_F^\times$, ${\rm N}_{E/F}(\beta) \in \mathcal{O}_F^\times$, and $E = F(\beta)$; if ${\rm char}(F) \neq 2$, then there is a $\beta \in E$ such that $\bar{\beta} = -\beta$ and $E = F(\beta)$. Fix such a $\beta$.

If ${\bf G} = {\rm U}_3$, the elements of the maximal torus $T = M$ are of the form $t = {\rm diag}(a,b,\bar{a}^{-1})$, $a \in E$, $b \in E^1 = \ker({\rm N}_{E/F})$. The elements of $U$ are of the form
\begin{equation*}
   n(x,z(x,y)) = \left( \begin{array}{ccc} 1 & x & z(x,y) \\ & 1 & \bar{x} \\ & & 1 \end{array} \right), \ x \in E, \, y \in F,
\end{equation*}
where $z(x,y) \in E$ satisfies the equation ${\rm Tr}_{E/F}(z(x,y)) = {\rm N}_{E/F}(x)$. In the case of a field extension $E/F$, this implies that
\begin{equation*}
   z(x,y) = \left\{ \begin{array}{ll}  y + \beta({\rm Tr}_{E/F}(\beta))^{-1}{\rm N}_{E/F}(x) & \text{ if } {\rm char}(F) = 2\\
   			     				\frac{1}{2} {\rm N}_{E/F}(x) + \beta y & \text{ if } {\rm char}(F) \neq 2
			      \end{array} \right. .
\end{equation*}
If $E = F \times F$, $x = (x_1,x_2)$, then
\begin{equation*}
   z(x,y) = (y,-y) + \beta x_1x_2.
\end{equation*}
To abbreviate, write $z = z(x,y)$. The non-trivial character $\psi$ of $F$ defines a character of $U$, also denoted by $\psi$, via the relationship
\begin{equation*}
   \psi(n(x,z)) = \overline{\psi} \circ {\rm Tr}_{E/F}(x).
\end{equation*}
Let $N_\alpha = \left\{ n(x,0) \vert x \in E \right\}$ and $N_{2\alpha} = \left\{ n(0,y) \vert y \in F \right\}$. Then $N = N_\alpha N_{2\alpha}$, $N_\alpha \cap N_{2\alpha} = \left\{ I_3 \right\}$. Notice that, as a topological group, $E = F \oplus \beta F$. Let $\mu$ be the Haar measure on $E$ given by
\begin{equation*}
   \int_{E} f(x_1 + \beta x_2) d\mu(x_1 + \beta x_2) = \int_{F} \int_{F} f(x_1 + \beta x_2) d\mu_\psi(x_1) d\mu_\psi(x_2).
\end{equation*}
Uniqueness of Haar measures gives
\begin{equation} \label{haarconst}
   \mu_{\psi_E} = c \mu.
\end{equation}
The measure $\mu_{\psi_E}$ of $E$ defines a Haar measure on $N_\alpha$, and the measure $\mu_\psi$ of $F$ defines a Haar measure on $N_{2\alpha}$. Then, fix the Haar measure on $N = N_\alpha N_{2\alpha}$ to satisfy $dn = c \, d\mu \, d\mu_\psi$.
If $E = F \times F$, $c=1$. Otherwise, $c$ depends on the level of $\psi_E$ with respect to $\psi$. If ${\rm char}(F) = 2$, the level of $\psi_E$ equals $l$ and $c=1$. It is an exercise to compute the constant $c$ in the case of ${\rm char}(F) \neq 2$, note that care must be taken for dyadic fields.

If ${\bf G} = {\rm SO}_3$, the elements of the maximal torus $T = M$ are of the form $t = {\rm diag}(a,1,a^{-1})$. The elements of $U$ are of the form
\begin{equation*}
   n(x) = \left( \begin{array}{ccc} 1 & \! -2x & \! -x^2 \\  & \, 1 & \, x \\  & & \, 1 \end{array} \right), \ x \in F.
\end{equation*}
This includes the case ${\rm char}(F) = 2$. The self-dual Haar measure $\mu_\psi$ of $F$ defines a Haar measure on $N = \left\{ n(x) \vert x \in F \right\}$. The character $\psi$ of $F$ is then extended to a character $\psi$ of $U$ satisfying $\psi(n(x)) = \overline{\psi}(x)$.

If ${\bf G} = {\rm GL}_2$, ${\rm SL}_2$, then the elements of $U$ are of the form
\begin{equation*}
   n(x) = \left( \begin{array}{cc} 1 & x \\ & 1 \end{array} \right), \ x \in E.
\end{equation*}
And, the character $\psi$ of $F$ is extended to a character of $U$ by setting $\psi(n(x)) = \psi_E(x)$. 

If ${\bf G} = {\rm U}_2$, then
\begin{equation*}
   n(x) = \left( \begin{array}{cc} 1 & x \\ & 1 \end{array} \right), \ x \in F.
\end{equation*}
Recall that $x$ is identified with $(x,x) \in E$ in the case $E = F \times F$. The character $\psi$ of $F$ is extended to a character of $U$ by setting $\psi(n(x)) = \psi(x)$.

Fix representatives of the Weyl group element $w_0$ appearing in the definition of the local coefficient:
\begin{equation*}
   w_0 = \left( \begin{array}{rr} 0 & 1 \\ -1 & 0 \end{array} \right), \text{ if } {\bf G} = {\rm GL}_2, \ {\rm SL}_2, \text{ or } {\rm U}_{2}, 
\end{equation*}
and
\begin{equation*}
   w_0 = \left( \begin{array}{ccc} 0 & \, 0 & 1 \\ 0 & \!\!\! -1 & 0 \\ 1 & \, 0 & 0 \end{array} \right), \text{ if } {\bf G} = {\rm SO}_{3}, \text{ or } {\rm U}_3.
\end{equation*}

\subsection{Compatibility with class field theory} 
\
\medskip

\noindent{\bf Proposition.} \emph{Let $\chi$, $\chi_1$, and $\chi_2$ be continuous characters of ${\rm GL}_1(E)$}.
\begin{itemize}
\item \emph{If ${\bf G} = {\rm GL}_2$ and $\pi$ is the smooth representation of $T$ given by $\pi({\rm diag}(t_1,t_2)) = \chi_1(t_1)\chi_2^{-1}(t_2)$, then
\begin{equation*}
   C_\psi(s,\pi,w_0) = \gamma(s,\chi_1 \chi_2,\psi).
\end{equation*}
}
\item \emph{If ${\bf G} = {\rm SO}_3$ and $\pi$ is the smooth representation of $T$ given by $\pi({\rm diag}(t,1,t^{-1})) = \chi(t)$, then
\begin{equation*}
   C_\psi(s,\pi,w_0) = \gamma(s,\chi^2,\psi).
\end{equation*}
}
\item \emph{If ${\bf G} = {\rm SL}_2$ and $\pi$ is the smooth representation of $T$ given by $\pi({\rm diag}(t,t^{-1})) = \chi(t)$, then
\begin{equation*}
   C_\psi(s,\pi,w_0) = \gamma(s,\chi,\psi).
\end{equation*}
}
\item \emph{If ${\bf G} = {\rm U}_2$ and $\pi$ is the smooth representation of $T$ given by $\pi({\rm diag}(t,\bar{t}^{-1})) = \chi(t)$, then
\begin{equation*}
   C_\psi(s,\pi,w_0) = \gamma(s,\chi \vert_{F^\times}, \psi).
\end{equation*}
}
\item \emph{If ${\bf G} = {\rm U}_3$ and $\nu$ is a continuous character of $E^1$, extended to a character of $E^\times$ via Hilbert's theorem~90. If $\pi$ is the smooth representation on $T$ given by $\pi({\rm diag}(t,z,\bar{t}^{-1})) = \chi(t) \nu(z)$, then
\begin{equation*}
   C_\psi(s,\pi,w_0) = \lambda(E/F,\overline{\psi}) \, \gamma_E(s,\chi \nu, \psi_E) \, \gamma(2s, \eta_{E/F} \chi \vert_{F^\times},\psi).
\end{equation*}
}
\end{itemize}

\noindent {\bf Remark.} \emph{When ${\bf G}$ is a quasi-split unitary group and $E = F \times F$, embed $F$ diagonally in $E$. The character $\chi$ is then of the form $\chi_1 \otimes \chi_2$, and $\chi \vert_{F^\times}$ is given by $\chi \vert_{F^\times}((x,x)) = \chi_1(x) \chi_2(x)$.
}

\medskip

\noindent{\it Proof.} The case of ${\rm U}_3$ is presented, with all details given in the case when $E = F \times F$ or $E/F$ is an extension of local fields with ${\rm char}(F) = 2$. With the above notation:
\begin{align*}
   \lambda_\psi(-s,w_0(\pi))({\rm A}(s,\pi,w_0)f_s^0)
   	& = \int_N {\rm A}(s,\pi,w_0)f_s^0(w_0^{-1}n(u,v)) \psi_E(u) dn(u,v)
\end{align*}
\begin{equation*}
   	= \int_N \int_N f_s^0(w_0^{-1}n(x,z)w_0^{-1}n(u,v)) \psi_E(u) dn(x,z) dn(u,v).
\end{equation*}
Using the Bruhat decomposition with $z \neq 0$ gives this equal to
\begin{align*}
	\int_N \int_{N - \left\{ I_3 \right\}} f_s^0 ( {\rm diag}(\bar{z}^{-1}, \bar{z}z^{-1},z) n(-x\bar{z}z^{-1},\bar{z}) & w_0 n(-xz^{-1},z^{-1}) n(u,v) ) \\
	& \times \psi_E(u) dn(x,z)  dn(u,v),
\end{align*}
and noticing that $n(-xz^{-1},z^{-1}) n(u,v) = n(u,v) n(-xz^{-1},z')$, for some $z' \in E$, makes this equal to
\begin{equation*}
	\int_N \left| \bar{z} \right|_E^{-(s+1)} \chi^{-1}(\bar{z}) \nu(z) \overline{\psi}_E(-xz^{-1}) dn(x,z) \left( \lambda_\psi(s,\pi)f_s^0 \right).
\end{equation*}
Therefore
\begin{equation*}
   C_\psi(s,\pi,w_0)^{-1} = \int_N \left| \bar{z} \right|_E^{-(s+1)} \chi^{-1}(\bar{z}) \nu(z) \psi_E(xz^{-1}) dn(x,z).
\end{equation*}
Consider the cases $E = F \times F$ or ${\rm char}(F) = 2$ (it is an exercise to adapt the following discussion to the remaining case when $E/F$ is an extension of local fields with ${\rm char}(F) \neq 2$). Notice that $z(x,(\beta + \bar{\beta})^{-1}x\bar{x}y) = (\beta + \bar{\beta})^{-1}x\bar{x}(y+\beta)$. Temporary notation: if $E = F \times F$ and $y \in F$, identify $y$ with $(y,-y) \in E$. Then, the integral
\begin{align*}
   \int_E \int_F & \left| \bar{z}(x,y) \right|_E^{-(s+1)} \chi^{-1}(\bar{z}(x,y)) \nu(z(x,y)) \psi_E(xz(x,y)^{-1}) d\mu_\psi(y) d\mu_{\psi_E}(x)
\end{align*}
is equal to
\begin{align*}
    \int_F \int_E \left| \overline{y+\beta} \right|_E^{-s-1} & \left| x \right|_E^{-2s-1} \chi^{-1}((\beta+\bar{\beta})^{-1} x \bar{x} (\overline{y+\beta})) \nu(y+\beta) \\
    					& \times \psi_E((\beta+\bar{\beta})(y+\beta)^{-1} \bar{x}^{-1}) d\mu_{\psi}(x) d\mu_{\psi_E}(y).
\end{align*}
Hence, after taking $x \mapsto (\beta+\bar{\beta})(\overline{y+\beta})^{-1}x$,
\begin{align*}
   C_\psi(s,\pi,w_0)^{-1} = \ & c \, \chi^{-1}(\beta+\bar{\beta}) \int_F \left| y+\beta \right|_E^{s-1} \chi(y+\beta) \nu(y+\beta) d\mu_\psi(y) \\
   & \times \int_E \left| x \right|_E^{-2s-1} \chi^{-1}({\rm N}_{E/F}(x)) \psi_E(\overline{x}^{-1}) d\mu_{\psi_E}(x),
\end{align*}
where $c$ is the constant of equation~\eqref{haarconst}. Now, notice that
\begin{align*}
   \chi^{-1} & (\beta+\bar{\beta}) \, \gamma(2s,\chi\vert_{F^\times},\psi)^{-1} \int_F \left| y+\beta \right|_E^{s-1} \chi(y+\beta) \nu(y+\beta) d\mu_\psi(y) \\
   & = \int_F \left| x \right|_F \int_F \left| yx + \beta x \right|_E^{s-1} \chi(yx+\beta x) \nu(yx+\beta x) \psi((\beta+\bar{\beta})x) d\mu_\psi(y) d\mu_\psi(x) \\
   & = \int_F \int_F \left| y+\beta x \right|_E^{s-1} \chi(y+\beta x) \nu(y + \beta x) \psi \circ {\rm Tr}_{E/F}(y+\beta x) d\mu_{\psi}(y) d\mu_{\psi}(x) \\
   & = c^{-1} \int_E \left| z \right|_E^{s-1} \chi(z) \nu(z) \psi_E(z) d\mu_{\psi_E}(z) = \gamma_E(s,\chi \nu,\psi_E)^{-1}.
\end{align*}
Also, taking $x \mapsto x^{-1}$ followed by $x \mapsto \bar{x}$, gives
\begin{align*}
   \int_E & \left| x \right|_E^{-2s-1} \chi^{-1}({\rm N}_{E/F}(x)) \psi_E(\bar{x}^{-1}) d\mu_{\psi_E}(x) = \gamma_E(2s,\chi \circ {\rm N}_{E/F},\psi_E)^{-1} \\
   & = \lambda(E/F,\psi) \, \gamma(2s,\chi\vert_{F^\times},\psi)^{-1} \gamma(2s,\eta_{E/F} \chi\vert_{F^\times},\psi)^{-1},
\end{align*}
where the Langlands $\lambda$-factor comes from equation~\eqref{langlandslambda}. Thus
\begin{equation*}
   C_\psi(s,\pi,w_0) = \lambda(E/F,\overline{\psi}) \gamma_E(s,\chi \nu,\psi_E) \gamma(2s,\eta_{E/F} \chi\vert_{F^\times},\psi).
\end{equation*}

\section{Principal series}

\subsection{Weyl group elements} For every simple root $\alpha$, fix a Weyl group element representative $w_\alpha$ in such a way that is in accordance with the semisimple rank one theory. Specifically, given a classical group ${\bf G} = {\bf G}_n$ ($n \geq 2$), for each $i$, $1 \leq i \leq n-1$, let $\alpha_i$ be the simple root $e_i - e_{i+1}$ in the Bourbaki notation. Embed ${\bf G}_{\alpha_i} = {\rm GL}_2$ in ${\bf G}$ so that
\begin{equation*}
   G_{\alpha_i} = \left\{ \left(	\begin{array}{ccccc}	I_{i-1} & & & & \\
   									 & g & & & \\
									 & & I_{2n-2i-2} & & \\
									 & & & \Phi_2 {}^t \bar{g}^{-1} \Phi_2 & \\
									 & & & & I_{i-1}
   					\end{array} \right) \vert \ g \in {\rm GL}_2(E) \right\} \subset G. 
\end{equation*}
Then, the representatives $w_{\alpha_i}$ are fixed for $1 \leq i \leq n-1$. The $w_\alpha$'s corresponding to the remaining simple roots are given by the embedding of~\eqref{G_1embed}, unless ${\bf G} = {\rm SO}_4$. In this latter case, let $w_{\alpha_1}$ correspond to the simple root $\alpha_1 = e_1 - e_2$ as before, and fix
\begin{equation*}
   w_{\alpha_2} = \left( \begin{array}{cccc} & & -1 & \\
   							  & & & 1 \\
							  1 & & & \\
							  & -1 & &
   			   \end{array} \right)
\end{equation*}
for the simple root $\alpha_2 = e_1 + e_2$.

Given an element $w$ of the Weyl group, it has a reduced decomposition into a product of simple roots. The representative of $w$ given by this product is independent of the reduced decomposition.

\subsection{Multiplicativity and Haar measures} For future reference, Langlands' lemma is recalled below for a quasi-split connected reductive group $\bf G$ defined over a field $k$. Shahidi's proof (see \S~2.1 of \cite{sha81}, originally written for non-archimedean local fields of characteristic zero), provides an algorithm to properly decompose the unipotent radical $\bf N$ of a parabolic subgroup $\bf P$ of $\bf G$ with respect to the decomposition of the corresponding Weyl group elements. If $\bf G$ is viewed as a reductive group scheme, the algorithm further extends to include quasi-split unitary groups with respect to a degree-$2$ finite \'etale algebra over $k$ (see 4.4.5 of \cite{conrad}).

Let $\bf G$ be a quasi-split connected reductive algebraic group defined over $k$. Let $\bf B$ be a borel subgroup of $\bf G$ with maximal split torus ${\bf A}_0$ and unipotent radical $\bf U$. Denote the roots of $\bf G$ with respect to ${\bf A}_0$ by $\Sigma$, the positive roots with respect to $\bf B$ by $\Sigma^+$, and the simple roots by $\Delta$. Parabolic subgroups $\bf P$ of $\bf G$ correspond to subsets $\theta$ of $\Delta$. Then ${\bf P} = {\bf P}_\theta$ has Levi decomposition ${\bf M}_\theta {\bf N}_\theta$. Let $P_\theta$, $M_\theta$, $N_\theta$ denote the groups of rational points. Let $W$ be the Weyl group of $\Sigma$, and for each $\alpha \in \Delta$ let $w_\alpha$ be its corresponding reflection. Given $\theta \subset \Delta$, let $W_\theta$ be the subgroup of $W$ generated by $w_\alpha$, $\alpha \in \theta$, and let $w_{l,\theta}$ be the longest element of $W_\theta$. Let $\Sigma_\theta$ be the set of roots spanned by $\theta$ and let $\Sigma_\theta^+ = \Sigma^+ \cap \Sigma_\theta$.

Two subsets $\theta$ and $\theta'$ of $\Delta$ are said to be associate if the set
\begin{equation*}
   W(\theta,\theta') = \left\{ w \in W \vert w(\theta) = \theta' \right\}
\end{equation*}
is non-empty. Given two such sets, the corresponding parabolic subgroups ${\bf P}_\theta$ and ${\bf P}_{\theta'}$ of $\bf G$ are associate. For an element $w \in W$, let $N_w = U \cap w N_{\theta}^-w^{-1}$ and let $\overline{N}_w = w^{-1}N_ww$.

\medskip

\noindent{\bf Lemma} \emph{Suppose $\theta$, $\theta' \subset \Delta$ are associate. Take $w \in W(\theta,\theta')$. Then, there exists a family of subsets $\theta_1, \ldots, \theta_n$ of $\Delta$, with $\theta_1 = \theta$ and $\theta_n = \theta'$, such that for every $i$, $1 \leq i \leq n-1$
\begin{itemize}
   \item[(i)] there exists a root $\alpha_i \in \Delta - \theta_i$ such that $\theta_{i+1}$ is the conjugate of $\theta_i$ in $\Omega_i = \theta_i \cup \left\{ \alpha_i \right\}$;
   \item[(ii)] if $w_i = w_{l,\Omega_i}w_{l,\theta_i}$ in $W(\theta_i,\theta_{i+1})$, then $w = w_{n-1} \cdots w_1$;
   \item[(iii)] if one sets $w_1' = w$ and $w_{i+1}' = w_i'w_i^{-1}$, then $w_n'=1$ and
      \begin{equation*}
         \overline{N}_{w_i'} = w_i^{-1} \overline{N}_{w_{i+1}'} w_i \rtimes \overline{N}_{w_i}.
      \end{equation*}
\end{itemize}
}

\bigskip

If $k$ is a non-archimedean local field, the semi-direct product of part~(iii) of the lemma is a semi-direct product of analytic $\mathfrak{p}$-adic groups. With the notation of the lemma, Haar measures are then normalized to satisfy
\begin{equation*}
   \int_{\overline{N}_{w_i'}} f(\bar{n}_i) \, d\bar{n}_i = \int_{\overline{N}_{w_{i+1}'} \times \overline{N}_{w_i}} f(w_i^{-1}\bar{n}_{i+1}'w_i\bar{n}_i) \, d\bar{n}_{i+1}' d\bar{n}_i. 
\end{equation*}
Now, let $\bf G$ is a classical group, $w = w_0$, and $\theta = \theta' = \Delta - \left\{ \alpha \right\}$ corresponds to the maximal Levi subgroup $\bf M$. Then, the lemma allows integrals over $N$ to be decomposed into a product of integrals over unipotent radicals corresponding to semisimple rank one groups. In this way, Haar measures are normalized so that they are in accordance with the semisimple rank one theory.

Let $X(M_\theta)$ be the set or rational characters of $M_\theta$ and let $\mathfrak{a}_{\theta,\mathbb{C}}^* = X(M_\theta) \otimes \mathbb{C}$. Given $w \in W(\theta,\theta')$, fix a decomposition $w = w_{n-1} \cdots w_1$ as in the lemma. Let $\nu_1 = \nu \in \mathfrak{a}_{\theta,\mathbb{C}}^*$, and let $\sigma_1 = \sigma$ be an irreducible generic representation of $M_\theta$. For $i$, $1 \leq i \leq n-1$, let $\nu_i = w_i(\nu_{i-1})$, $\sigma_i = w_{i-1}(\sigma_{i-1})$. Then, multiplicativity for the local coefficient can be stated as follows:
\begin{equation*}
   C_\psi(\nu,\sigma,w) = \prod_{i=1}^{n-1} C_\psi(\nu_i,\sigma_i,w_i).
\end{equation*}

\subsection{Principal series} For every integer $m \geq 1$, let ${\bf B}_m$ be the Borel subgroup of ${\rm GL}_m$, whose group $B_m$ of $F$-rational points consists of upper triangular matrices. Let $\chi_i$, $\mu_j$, $1 \leq i \leq n_1$, $1 \leq j \leq n_2$, be continuous characters of $F^\times$. Let $\pi_1$ be a $\psi$-generic constituent of
\begin{equation*}
   {\rm ind}_{B_{n_1}}^{{\rm GL}_{n_1}(F)}(\chi_1 \otimes \cdots \otimes \chi_{n_1}),
\end{equation*}
and $\pi_2$ a $\overline{\psi}$-generic constituent of
\begin{equation*}
   {\rm ind}_{B_{n_2}}^{{\rm GL}_{n_2}(F)}(\mu_1 \otimes \cdots \otimes \mu_{n_2}).
\end{equation*}
The algorithm of the previous section gives a decomposition of the Weyl group element $w_0 = w_{n_1n_2} w_{n_1n_2-1} \cdots w_1$, where the $w_i's$, $1 \leq i \leq n_1n_2$, correspond to positive simple roots. Notice that if a smooth representation $\tau$ of ${\rm GL}_m(F)$ is $\psi$-generic, then $\tilde{\tau}$ is $\overline{\psi}$-generic. The multiplicativity property of the local coefficient for the $\psi$-generic representation $\pi = \pi_1 \otimes \tilde{\pi}_2$ of $M$ reads as follows
\begin{equation*}
   C_\psi(s,\pi,w_0) = \prod_{i,j} C_\psi(s,\chi_i \otimes \mu_j^{-1}, w_0').
\end{equation*}
Here $C_{\psi}(s,\chi_i \otimes \mu_j^{-1},w_0')$ is a local coefficient of ${\rm GL}_2$ and $w_0'$ is a Weyl group element representative corresponding to the element denoted $w_0$ in the rank one case. Furthermore, there is a connection to the local theory of Tate's thesis due to Proposition~3.2.

\medskip

\noindent{\bf 4.4. Proposition.} \emph{Let $\pi_1$ be a $\psi$-generic constituent of
\begin{equation*}
   {\rm ind}_{B_{n_1}}^{{\rm GL}_{n_2}(F)}(\chi_1 \otimes \cdots \otimes \chi_{n_1}),
\end{equation*}
and let $\pi_2$ be a $\overline{\psi}$-generic constituent of
\begin{equation*}
   {\rm ind}_{B_{n_2}}^{{\rm GL}_{n_2}(F)}(\mu_1 \otimes \cdots \otimes \mu_{n_2}).
\end{equation*}
Then, $\pi = \pi_1 \otimes \tilde{\pi}_2$ is a $\psi$-generic representation of $M \subset {\rm GL}_{n_1+n_2}(F)$ and the corresponding local coefficient satisfies
\begin{equation*}
   C_\psi(s,\pi_1 \otimes \tilde{\pi}_2,w_0) = \prod_{i,j} \gamma(s,\chi_i \mu_j,\psi).
\end{equation*}
This product is denoted by $\gamma(s,\pi_1 \times \pi_2,\psi)$, a Rankin-Selberg $\gamma$-factor.
}

\medskip

Similarly, the algorithm applied to local coefficients of Siegel Levi subgroups of the classical groups and principal series for ${\rm GL}_n$ gives the following proposition. Notice that Proposition~4.4 can be extended to include $E = F \times F$, in particular, the semisimple rank one case of ${\rm GL}_2$ is needed to prove Proposition~4.5 for unitary groups in this case.

\medskip

\noindent{\bf 4.5. Proposition} \emph{When $E$ is a non-archimedean local fied, let $\pi$ be a $\psi$-generic constituent of
\begin{equation*}
   {\rm ind}_{B_n}^{{\rm GL}_n(E)}(\chi_1 \otimes \cdots \otimes \chi_{n}).
\end{equation*}
In the case of the separable quadratic algebra $E = F \times F$, let $\pi = \pi_1 \otimes \pi_2$, where $\pi_1$ is a $\psi$-generic constituent of
\begin{equation*}
   {\rm ind}_{B_n}^{{\rm GL}_n(F)}(\chi_{1,1} \otimes \cdots \otimes \chi_{1,n}),
\end{equation*}
and $\pi_2$ is a $\psi$-generic constituent of
\begin{equation*}
   {\rm ind}_{B_n}^{{\rm GL}_n(F)}(\chi_{2,1} \otimes \cdots \otimes \chi_{2,n}).
\end{equation*}
In the latter case, write $\chi_i = \chi_{1,i} \otimes \chi_{2,i}$, for each $i$, $1 \leq i \leq n$. Then $\pi$ is an irreducible $\psi_{M}$-generic representation of $M$ and the local coefficient satisfies
\begin{equation*}
   C_\psi(s,\pi,w_0) = \left\{ \begin{array}{ll} \prod_{i=1}^n \gamma(s,\chi_i^2,\psi) \prod_{i<j} \gamma(s,\chi_i \chi_j,\psi) & \text{if } {\bf G} = {\rm SO}_{2n+1} \\
   								   \prod_{i=1}^n \gamma(s,\chi_i,\psi) \prod_{i<j} \gamma(2s,\chi_i \chi_j,\psi) & \text{if } {\bf G} = {\rm Sp}_{2n} \\
								   \prod_{i<j} \gamma(s,\chi_i \chi_j,\psi) & \text{if } {\bf G} = {\rm SO}_{2n} \\
								   \prod_{i=1}^n \gamma(s,\chi_i\vert_{F^\times},\psi) \prod_{i<j} \gamma_E(s,\chi_i \chi_j^{{\rm conj.}},\psi_E) 
								   & \text{if } {\bf G} = {\rm U}_{2n}
					\end{array} \right. .
\end{equation*}
}

\emph{Let $\nu'$ be a smooth representation of $E^1$, and extend it to a continuous character $\nu$ of $E^\times$ by setting $\nu = \nu'(\bar{z}z^{-1})$. If ${\bf G} = {\rm U}_{2n+1}$, then $\pi \otimes \nu$ is a $\psi_M$-generic representation of $M$ and the local coefficient satisfies
\begin{align*}
   C_\psi(s,\pi \otimes \nu,w_0) = \lambda(E/F,\overline{\psi})^{n} & \prod_{i=1}^n \gamma_E(s,\chi_i\nu,\psi_E) \gamma(2s,\eta_{E/F}\chi_i\vert_{F^\times},\psi)  \\
   												& \times \prod_{i<j} \gamma_E(2s,\chi_i \chi_j^{\rm conj.},\psi_E). 
\end{align*}
}

\noindent \emph{Notation:} if $\bf H$ is a group defined over $E$, and $\tau$ is a representation of the group of $E$-rational points $H$ of $\bf H$, then $\tau^{\rm conj.}$ denotes the representation of $H$ defined via conjugation, i.e., $\tau^{\rm conj.}(g) = \tau(\bar{g})$.

Assume $E = F \times F$. Notice that, if ${\bf G} = {\rm U}_{2n}$, then
\begin{equation*}
   C_\psi(s,\pi_1 \otimes \pi_2,w_0) = \prod_{i,j} \gamma(s,\chi_{1,i} \chi_{2,j},\psi) = \gamma(s, \pi_1 \times \pi_2, \psi).
\end{equation*}
If ${\bf G} = {\rm U}_{2n+1}$, then $\nu$ is given by $\nu((x,y)) = \nu_0^{-1}(x) \nu_0(y)$, with $\nu_0: F^\times \rightarrow \mathbb{C}$ continuous, and
\begin{align*}
   C_\psi(s,\pi,w_0)	& = \prod_{i=1}^n \gamma_E(s,\chi_i\nu,\psi_E) \prod_{i,j} \gamma(2s,\chi_{1,i} \chi_{2,j},\psi) \\
   				& = \prod_{i=1}^n \gamma(s,\chi_{1,i} \nu_0^{-1},\psi) \gamma(s,\chi_{2,i}\nu_0,\psi) \prod_{i,j} \gamma(2s,\chi_{1,i} \chi_{2,j},\psi) \\
				& = \gamma(s,\pi_1 \times \nu_0^{-1},\psi) \gamma(s,\pi_2 \times \nu_0,\psi) \gamma(2s,\pi_1 \times \pi_2,\psi).
\end{align*}

\section{Global functional equation of local coefficients}

Let $k = K$ be a global function field if ${\bf G}$ is split, and let $K/k$ be a separable quadratic extension of global function fields if ${\bf G}$ is quasi-split. Let $\mathbb{A}_k$, respectiverly $\mathbb{A}_K$, denote the ring of ad\`eles of $k$, respectively $K$. Let $\psi = \otimes_v \psi_v$ be a non-trivial continuous character of $k \backslash \mathbb{A}_k$ such that each $\psi_v$ is unramified. In the case of a quadratic extension $K/k$ and a place $v$ that splits in $K$, let $K_v$ denote the separable quadratic algebra $k_v \times k_v$.
If ${\bf G} = {\bf G}_m$, $m=2n$ or $m=2n+1$, is a unitary group and $v$ is any place of $k$, let ${\bf G}_v$ be the unitary group over $K_v/k_v$ whose group of rational points is
\begin{equation*}
G_v = \left\{ g \in {\rm GL}_m({K_v}) \vert h_m(gx) = h_m(x) \text{ for every } x \in K_v^m  \right\}.
\end{equation*}
Let $\mathcal{K} = \prod_v \mathcal{K}_{v}$ be a maximal compact subgroup of ${\bf G}(\mathbb{A}_k)$. For every $v$, assume $\mathcal{K}_v$ is a special maximal compact subgroup of $G_v$, hyperspecial for almost all $v$. Then ${\bf G}(\mathbb{A}_k) = \otimes' G_v$, where the restricted direct product is relative to the subgroups $\mathcal{K}_{v}$.

Proceeding as in the local theory, the character $\psi$ is then extended to a character $\psi_M$ of ${\bf U}_M(\mathbb{A}_k)$, where ${\bf U}_M = {\bf M} \cap {\bf U}$, trivial on ${\bf U}_M(k)$. Then $\psi_M$, in turn, is extended to a character $\psi$ of ${\bf U}(\mathbb{A}_k)$, trivial on ${\bf U}(k)$. The local coefficients satisfy a functional equation, proved using the theory of Eisenstein series for globally generic automorphic representations. Notice that a system of Weyl group element representatives fixed in ${\bf G}(k)$, fixes a system of Weyl group elements at every place of $k$. We fix such a global system, which is in accordance with the local theory.

\medskip

\subsection{Theorem}(The crude functional equation). \emph{Let $\pi = \otimes' \, \pi_w$ be a globally $\psi_M$-generic cuspidal automorphic representation of ${\bf M}(\mathbb{A}_k) \subset {\bf G}(\mathbb{A}_k)$. If $\bf G$ is a quasi-split unitary group and the place $v$ in $k$ splits, let $\pi_v = \pi_{w_1} \otimes \pi_{w_2}$, where $w_1$, $w_2$ are the places of $K$ lying above $v$, otherwise write $\pi_w = \pi_v$. Let $S$ be a finite set of places such that $\pi_v$ is unramified for $v \notin S$. The local coefficients satisfy the functional equation:
\begin{equation*}
   \prod_{i=1}^{m_r}L^S(is,\pi,r_i) = \prod_{v \in S} C_{\psi_v}(s,\pi_v,w_0) \prod_{i=1}^{m_r} L^S(1-is,\tilde{\pi}_v,r_i).
\end{equation*}
}

\subsection{Partial $L$-functions} \label{partial}

To explain the notation in the crude functional equation, let ${}^LG$ be the $L$-group of ${\bf G}$ over $k$. Given an irreducible admissible automorphic representation $\pi = \otimes' \pi_v$ of ${\bf G}(\mathbb{A}_k)$, $\pi_v$ is unramified at almost every $v$. Its contragredient representation is denoted by $\tilde{\pi}$. 

Let $\rho$ be a representation of ${}^LG$, i.e., a continuous homomorphism $\rho : {}^LG \rightarrow {\rm GL}_m(\mathbb{C})$ whose restriction to ${}^LG^0$ is a morphism of complex Lie groups. Let $S$ be a finite set of places of $k$ such that $\pi_v$ is unramified for $v \notin S$. For each $v \notin S$, the Satake isomorphism determines a semisimple class $\hat{A}_v$ in ${}^LG$. The local $L$-function attached to $\pi_v$ and $\rho_v$, $v \notin S$, is defined by
\begin{equation*}
   L(s,\pi_v,\rho_v) = {\rm det}(I - \rho_v(\hat{A}_v)q_v^{-s})^{-1}.
\end{equation*}
Then, the partial $L$-function is defined by
\begin{equation*}
   L^S(s,\pi,\rho) = \prod_{v \notin S} L(s,\pi_v,\rho_v).
\end{equation*}
It is absolutely convergent in some right half plane (see \S~13 of \cite{borel}).

Given a positive integer $m$, let $\rho_m$ denote the standard representation of ${\rm GL}_m(\mathbb{C})$. Label $r = r_1$ to be $\rho_{n_1} \otimes \tilde{\rho}_{n_2}$, ${\rm Sym}^2 \rho_n$, or $\wedge^2 \rho_n$ depending on wether ${\bf G} = {\rm GL}_{n_1+n_2}$, ${\rm SO}_{2n+1}$, or ${\rm SO}_{2n}$. If ${\bf G} = {\rm U}_{2n}$, notice that the $L$-group is given by ${}^L{\rm Res}_{K/k}{\rm GL}_n = {\rm GL}_n(\mathbb{C}) \times {\rm GL}_n(\mathbb{C}) \rtimes \mathcal{W}_{k}$, where the Weil group acts via the Galois group ${\rm Gal}(K/k) = \left\{ 1,\theta \right\}$. In this case, let $r = r_1 = r_{\mathcal{A}}$ be the Asai representation from ${}^L{\rm Res}_{K/k}{\rm GL}_n$ to ${\rm GL}_{n^2}(\mathbb{C})$, given by
\begin{equation*}
   r_\mathcal{A}(x,y,1) = x \otimes y, \text{ and } r_\mathcal{A}(x,y,\theta) = y \otimes x.
\end{equation*}
Also, if ${\bf G} = {\rm Sp}_2 = {\rm SL}_2$, let $r = r_1 = \rho_1$. In all of these cases $r $ is irreducible and $m_r = 1$ in Theorem~5.1. Proposition~4.4, together with equation~\eqref{1unramified}, and the Satake parametrization applied to $\pi_v$, $v \notin S$, yields
\begin{equation*}
   C_{\psi_v}(s,\pi_v,w_0) L(s,\pi_v,r_v) = L(1-s,\tilde{\pi}_v,r_v).
\end{equation*}
In the case ${\bf G} = {\rm SO}_{2n}$, all factors in the previous equation are interpreted to be trivial when $n=1$.

Now, if ${\bf G} = {\rm Sp}_{2n}$ with $n \geq 2$, let $r_1 = \rho_n$ and $r_2 = \wedge^2 \rho_n$. Then, for $v \notin S$,
\begin{equation*}
   C_{\psi_v}(s,\pi_v,w_0) L(s,\pi_v,r_1) L(2s,\pi_v,r_2) = L(1-s,\tilde{\pi}_v,r_1) L(1-2s,\tilde{\pi}_v,r_2).
\end{equation*}
If ${\bf G} = {\rm U}_{2n+1}$, let $r_1 = \rho_n \otimes \tilde{\rho}_1$ and $r_2 = r_\mathcal{A}$. Then, for $v \notin S$,
\begin{align*}
   \lambda(K_v/k_v,\pi_v)^n C_\psi(s,\pi_v,w_0) & L(s,\pi_v \otimes \nu_v^{-1},r_{1,v}) L(2s,\pi_v \otimes \eta_{K_v/k_v},r_{2,v}) \\
   				= & L(1-s,\tilde{\pi}_v \otimes \nu_v,r_{1,v}) L(1-2s,\tilde{\pi}_v \otimes \overline{\eta}_{K_v/k_v},r_{2,v}).
\end{align*}
In these two latter cases let $r = r_1 \oplus r_2$ and $m_r = 2$ in Theorem~5.1.

\subsection{Proof of the crude functional equation}

We now extend the discussion of \S~5 of \cite{lomeli2009}, which is written for split groups, to include the cases at hand concerning the quasi-split unitary groups. Let $\pi$ be a globally $\psi$-generic cuspidal representation of ${\bf M}(\mathbb{A}_k)$. The irreducible constituents of the globally induced representation of $G$ given by the restricted direct product
\begin{equation*}
   {\rm I}(s,\pi) = \otimes' {\rm I}(s,\pi_v),
\end{equation*}
are automorphic representations $\Pi = \otimes' \Pi_v$ of $G$ such that the representation $\pi_v$ has $\mathcal{K}_v$-fixed vectors for almost all $v$. The space of ${\rm I}(s,\pi)$ is denoted by ${\rm V}(s,\pi)$. The restricted tensor product is taken with respect to functions $f_{v,s}^0$ that are fixed under the action of $\mathcal{K}_v$.

Since $\pi$ is globally $\psi$-generic, by definition, there is a cusp form $\varphi$ in the space of $\pi$ such that
\begin{equation*}
   W_{M,\varphi}(m) = \int_{{\bf U}_M(K) \backslash {\bf U}_M(\mathbb{A}_k)} \varphi(um) \overline{\psi}(u) \, du \neq 0.
\end{equation*}
It is possible to extend $\varphi$ to a cusp form $\tilde{\varphi}$ defined on ${\bf U}(\mathbb{A}_k){\bf M}(K) \backslash {\bf G}(\mathbb{A}_k)$, and define a function $\Phi_s$ through $\tilde{\varphi}$ in such a way that the Eisenstein series
\begin{equation*}
   E(s,\Phi_s,g) = \sum_{\gamma \in {\bf P}(k) \backslash {\bf G}(k)} \Phi_s(\gamma g)
\end{equation*}
satisfies
\begin{equation} \label{globaleswhittaker}
   E_\psi(s,\Phi_s,g,P) = \prod_v \lambda_{\psi_v}(s,\pi_v)({\rm I}(s,\pi_v)(g_v)f_{s,v}),
\end{equation}
with $f_s \in {\rm V}(s,\pi)$. Here $E_\psi(s,\Phi_s,g)$ denotes the Fourier coefficient
\begin{equation*}
   E_\psi(s,\Phi_s,g) = \int_{{\bf U}(K) \backslash {\bf U}(\mathbb{A}_k)} E(s,\Phi_s,ug) \overline{\psi}(u) \, du.
\end{equation*}
The global intertwining operator ${\rm M}(s,\pi)$ is defined by
\begin{equation*}
   {\rm M}(s,\pi,w_0)f(g) = \int_{{\bf N'}(\mathbb{A}_k)} f(w_0^{-1}ng) \, dn, 
\end{equation*}
where $f \in {\rm V}(s,\pi)$ and ${\bf N'}$ is the unipotent radical of the standard parabolic ${\bf P}'$ with Levi ${\bf M}' = w_0 {\bf M} w_0^{-1}$. It is the product of local intertwining operators
\begin{equation*}
   {\rm M}(s,\pi) = \prod_v {\rm A}(s,\pi_v,w_0).
\end{equation*}
It is a meromorphic operator, which is rational on $q^{- s}$ (Proposition~IV.1.12 of \cite{mw}). Property~\eqref{globaleswhittaker} gives
\begin{align*} \label{globaleswhittaker2}
   E_\psi&(-s,{\rm M}(s,\pi)\Phi_s,g,P') \\
   & = \prod_v \lambda_{\psi_v}(-s,w_0(\pi_v))({\rm I}(-s,w_0(\pi_v))(g_v) {\rm A}(s,\pi_v,w_0)f_{s,v}).
\end{align*}
It is known that the above Eisenstein series, and its Fourier coefficients, are rational functions on $q^{-s}$. The argument for proving this fact is due to Harder \cite{harder}.

Fourier coefficients of Eisenstein series satisfy the functional equation:
\begin{equation*}
   E_\psi(-s,{\rm M}(s,\pi)\Phi_s,g,P') = E_\psi(s,\Phi_s,g,P).
\end{equation*}
And, equation \eqref{globaleswhittaker} gives
\begin{align*}
   E_\psi(s,\Phi_s,e,P) & = \prod_v \lambda_{\psi_v}(s,\pi_v)f_{s,v} \\
   E_\psi(-s,{\rm M}(s,\pi)\Phi_s,e,P') & = \prod_v \lambda_{\psi_v}(-s,w_0(\pi_v)){\rm A}(s,\pi_v,w_0)f_{s,v}.
\end{align*}
Then, the Casselman-Shalika formula for unramified quasi-split groups, Theorem~5.4 of \cite{cs}, allows one to compute the Whittaker functional when $\pi_v$ is unramified:
\begin{equation*} \label{csformula}
   \lambda_{\psi_v}(s,\pi_v) f_{s,v}^0 = \prod_{i=1}^{m_r} L(1+is,\pi_v,r_i)^{-1} f_{s,v}^0(e_v).
\end{equation*}
Also, for $v \notin S$, the intertwining operator gives a function ${\rm A}(s,\pi_v,w_0) f_{s,v}^0 \in {\rm I}(-s,w_0(\pi_v))$ satisfying
\begin{equation*} \label{interformula}
   {\rm A}(s,\pi_v,w_0) f_{s,v}^0(e_v) = \prod_{i=1}^{m_r} \dfrac{L(is,\pi_v,r_i)}{L(1+is,\pi_v,r_i)} f_{s,v}^0(e_v).
\end{equation*}
This equation is established by means of the multiplicative property of the intertwining operator, which reduces the problem to semisimple rank one cases.

Finally, combining the last five equations together gives
\begin{equation*}
   \prod_{i=1}^{m_r} L^S(is,\pi,r_i) = \prod_{v \in S} \dfrac{\lambda_{\psi_v}(s,\pi_v)f_{s,v}}{\lambda_{\psi_v}(-s,w_0(\pi_v)){\rm A}(s,\pi_v,w_0)f_{s,v}} \prod_{i=1}^{m_r} L^S(1-is,\tilde{\pi},r_i).
\end{equation*}
Which establishes the crude functional equation of Theorem~5.1.

\subsection{Corollary}(Rationality of partial $L$-functions). \emph{For each $i$, $1 \leq i \leq m_r$, $L^S(s,\pi,r_i)$ has a meromorphic continuations to rational functions on $q^{-s}$.}

\medskip

\noindent\emph{Proof.} This now follows from the proof of the crude functional equation together with the result of Harder on the rationality of $E_\psi(s,\Phi_s,g)$ (see the paragraph preceding Theorem~1.1.6 of \cite{harder}).

\section{On certain $\gamma$-factors} \label{gammafactors}

A refinement of the crude functional equation satisfied by local coefficients leads to a definition of $\gamma$-factors. Although the results of the previous sections are for $\psi$-generic representations, it is possible to define local coefficients for generic representations in general. This is useful, since it leads to an explicit formula explaining the behavior of $\gamma$-factors as the non-trivial character $\psi$ varies.

Thanks to joint work with G. Henniart \cite{hl11,hlarxiv}, it is possible to uniquely characterize $\gamma$-factors in terms of local properties and their connection to the global theory via the functional equation. In positive characteristic, $\gamma$-factors and related local factors agree, via the local Langlands correspondence, with those defined by Deligne and Langlands \cite{deligne,langlands}.

\subsection{Notation.} Recall that, if $\bf G$ is a classical group, then $\bf M$ is taken to be a Siegel Levi subgroup of $\bf G$. Let $\mathscr{L}({\bf M},{\bf G})$ be the class whose objects are triples $(E/F,\pi,\psi)$ consisting of:
\begin{itemize}
   \item a non-archimedean local field $F$;
   \item $E = F$ if ${\bf G}$ is split, $E$ a degree-2 finite \'etale algebra over $F$ if ${\bf G}$ is a quasi-split unitary group;
   \item an irreducible generic representation $\pi$ of $M$;
   \item a non-trivial continuous character $\psi$ of $F$.
\end{itemize}
Let $\mathscr{L}(p,{\bf M},{\bf G})$ be the class consisting of $(E/F,\pi,\psi) \in \mathscr{L}({\bf M},{\bf G})$ with $F$ of characteristic $p$.

When $\bf G$ is split, it is convenient to simply write $(F,\pi,\psi)$ for an object of $\mathscr{L}({\bf M},{\bf G})$. In the case of a general linear group, ${\bf M} = {\rm GL}_{n_1} \times {\rm GL}_{n_2}$. The degree of a triple $(F,\pi,\psi) \in \mathscr{L}({\rm GL}_{n_1} \times {\rm GL}_{n_2},{\rm GL}_{n_1+n_2})$ is defined to be $(n_1,n_2)$. To ease the notation in this case, write $\mathscr{L}$ instead of $\mathscr{L}({\rm GL}_{n_1} \times {\rm GL}_{n_2},{\rm GL}_{n_1+n_2})$. Whenever $(F,\pi,\psi) \in \mathscr{L}$ is of degree $(n_1,n_2)$, it is assumed that the representation $\pi$ of $M$ is of the form $\pi = \pi_1 \otimes \tilde{\pi}_2$, with generic representations $\pi_1$ of ${\rm GL}_{n_1}(F)$ and $\pi_2$ of ${\rm GL}_{n_2}(F)$.

If $\bf G$ is a classical group, the Siegel Levi is either ${\rm GL}_n$ or ${\rm GL}_n \times {\rm U}_1$. In these cases, call $n$ the degree of a triple $(E/F,\pi,\psi) \in \mathscr{L}({\bf M},{\bf G})$. In the case of a triple $(E/F,\pi,\psi) \in \mathscr{L}({\bf M},{\rm U}_{2n+1})$, the representation $\pi$ is of the form $\pi = \pi' \otimes \nu'$, where $\nu'$ is a smooth representation of ${\rm U}_1(E)$.

\subsection{} Let $(E/F,\pi,\psi) \in \mathscr{L}({\bf M},{\bf G})$, where $\pi$ is a $\chi_M$-generic representation of $M$ with Whittaker functional $\lambda_M$. The non-degenerate character $\chi_M$ of $U_M$ is extended to a non-degenerate character $\chi$ of $U$ so that they are $w_0$-compatible, i.e.,
\begin{equation*}
   \chi_M(n) = \chi(w_0 n w_0^{-1}), n \in U_M.
\end{equation*} 
Notice that $\psi$ and $\psi_M$ are $w_0$-compatible. The local coefficient for $\chi_M$-generic representations is defined via the uniqueness property of Whittaker functionals and the relationship
\begin{equation*}
   \lambda_{\chi}(s,\pi) = C_{\chi}(s,\pi,w_0) \lambda_\chi(-s,w_0(\pi)) {\rm A}(s,\pi,w_0),
\end{equation*}
where
\begin{equation*}
   \lambda_\chi(s,\pi)f_s = \int_{N'} \lambda_{\chi_M}(f_s(w_0^{-1}n))\overline{\chi}(n) dn, \quad f_s \in {\rm I}(s,\pi).
\end{equation*}

Now, there is a $t \in {\bf T}(\overline{F})$ such that $\chi \circ {\rm Ad}(t)(n) = \chi(t^{-1} n t) = \psi(n)$ , $n \in N$; ${\rm Ad}(t)$ is $F$-rational (see the discussion in pp. 282-283 of \cite{sha90}). If $f \in {\rm I}(s,\pi)$, then $f_t \in {\rm I}(s,\pi_t)$, where $\pi_t(g) = \pi(t^{-1}gt)$, $f_t(g) = f(t^{-1}gt)$. The representation $\pi_t$ is then $\psi$-generic and belongs to the equivalence class of $\pi$.

Assume now that $\pi$ is $\psi^a$-generic, with $a \in F^\times$. In all the cases at hand it is possible to produce a $t$ such that $a_0^{-1} = w_0^{-1}t^{-1}w_0t$ is an element of the center of ${\bf M}(F)$. It is possible to obtain an explicit formula explaining the behavior of the local coefficient when it is defined via different non-trivial characters $\psi^a$, $a \in F^\times$. To see this, let $a_t \in \mathbb{C}^\times$, be the module of the automorphism $n \mapsto tnt^{-1}$. Then
\begin{equation} \label{psidependence}
   C_\psi(s,\pi_t,w_0) = a_t \, \omega_\pi(a_0) \left| {\rm det}(a_0) \right|_F^{is} \delta_P^{\frac{1}{2}}(a_0) C_\chi(s,\pi,w_0),
\end{equation}
where $i = 1$ or $1/2$ depending on equation~\eqref{eq:ind}, and $\omega_\pi$ the central character of $\pi$. Define
\begin{equation} \label{Cprime}
   C_\psi'(s,\pi_t,w_0) = a_t^{-1} \lambda(E/F,w_0)^{-1} C_\psi(s,\pi_t,w_0),
\end{equation}
where $\lambda(E/F,w_0) = \lambda(E/F,\overline{\psi})^n$ if ${\bf G} = {\rm U}_{2n+1}$, and it is equal to $1$ otherwise.

\subsection{Rankin-Selberg $\gamma$-factors} Let $(F,\pi,\psi) \in \mathscr{L}$, where $\pi = \pi_1 \otimes \tilde{\pi}_2$ is of degree $(n_1,n_2)$. Take $t$ as above such that $\pi_t$ is $\psi$-generic. Define
\begin{equation*}
   \gamma(s,\pi_1 \otimes \tilde{\pi}_2, r, \psi) = C'_\psi(s,\pi_t,w_0).
\end{equation*}
To simplify notation in this case, write $\gamma(s,\pi_1 \otimes \tilde{\pi}_2,\psi)$ for $\gamma(s,\pi_1 \otimes \tilde{\pi}_2,r,\psi)$. For the smooth generic representations $\pi_1$ of ${\rm GL}_{n_1}(F)$ and $\pi_2$ of ${\rm GL}_{n_2}(F)$, local Rankin-Selberg $\gamma$-factors are defined in \cite{jpss}:
\begin{equation*}
   \gamma(s,\pi_1 \times \pi_2,\psi) = \varepsilon(s,\pi_1 \times \pi_2,\psi) \dfrac{L(1-s,\tilde{\pi}_1 \times \tilde{\pi}_2)}{L(s,\pi_1 \times \pi_2)}.
\end{equation*}
The main result of \cite{sha84} establishes the equality of $\gamma$-factors
\begin{equation} \label{ls-rs}
   \gamma(s,\pi_1 \otimes \tilde{\pi}_2,\psi) = \gamma(s,\pi_1 \times \pi_2,\psi)
\end{equation}
using local methods. Notice that the choice of Weyl group elements of \S~4 gives
\begin{equation*}
   w_0 = w_l w_{l,M} = \left( \begin{array}{cc} & (-1)^{n_2+1} I_{n_2} \\ (-1)^{n_1+n_2+1}I_{n_1} \end{array} \right).
\end{equation*}
It is also possible to give a different proof of \eqref{ls-rs} via the local-to-global result of \cite{hl11} by giving a characterization of Rankin-Selberg $\gamma$-factors. An interesting consequence of this equality is the following stability property of $\gamma$-factors, known for Rankin-Selberg products (see \cite{js}):

\medskip

\noindent{\bf Proposition}. \emph{Let $F$ be a non-archimedean local field. Let $\pi_1$ and $\pi_2$ be irreducible generic representations of ${\rm GL}_{n_1}(F)$ with $\omega_{\pi_1} = \omega_{\pi_2}$, and let $\tau$ be an irreducible generic representation of ${\rm GL}_{n_2}(F)$. There is an integer $A$ with the following property: If $\chi$ is a character of $F^\times$ with conductor $\mathfrak{p}^a$, $a \geq A$, then
\begin{equation*}
   \gamma(s,(\pi_1 \otimes \chi) \otimes \tilde{\tau},\psi) = \gamma(s,(\pi_2 \otimes \chi) \otimes \tilde{\tau},\psi).
\end{equation*}
}
\medskip

The Godement-Jacquet $\gamma$-factor corresponding to a principal $L$-function of \cite{j}, for an irreducible generic representation $\pi$ of ${\rm GL}_n(F)$, can be obtained by setting
\begin{equation*}
   \gamma(s,\pi,\psi) = \gamma(s,\pi \otimes \mathds{1},\psi),
\end{equation*}
where $\mathds{1}$ is the trivial representation on ${\rm GL}_1(F)$. The above proposition gives stability of $\gamma$-factors for representations $\pi_1$ and $\pi_2$ of ${\rm GL}_n(F)$ sharing the same central character:
\begin{equation*}
   \gamma(s,\pi_1 \otimes \chi,\psi) = \gamma(s,\pi_2 \otimes \chi,\psi),
\end{equation*}
where $\chi$ is a highly ramified character.
\subsection{Definition/induction for the classical groups} Proposition~4.5, equation~\eqref{1unramified}, and the discussion in \S~5.2 show that $\gamma$-factors are compatible with the definition using $L$-groups via the Satake parametrization for unramified principal series. Given an irreducible $\chi_M$-generic representation of $M$, take $t$ so that $\chi \circ {\rm Ad}(t) = \psi$ as in section \S~6.1. With the notation of \S~\ref{partial}, whenever $r = r_1$ is irreducible, define
\begin{equation*}
   \gamma(s,\pi,r,\psi) = C'_\psi(s,\pi_t,w_0).
\end{equation*}
This is the main induction step of the Langlands-Shahidi method for the classical groups. Since both local and global Langlands conjectures are known in the case of ${\rm GL}_n$ in positive characteristic, it is possible to obtain an equality with the corresponding representations on the galois side.

\medskip

\noindent{\bf Theorem.} \emph{Let ${\bf G} = {\rm SO}_{2n+1}$, ${\rm SO}_{2n}$ or ${\rm U}_{2n}$, and let $(E/F,\pi,\psi) \in \mathscr{L}({\bf M},{\bf G})$ be of degree $n$, with $F$ of positive characteristic. Let $\sigma$ be the $n$-dimensional Frobenius semisimple Weil-Deligne representation of the Weil group $\mathcal{W}_E$ corresponding to $\pi$ under the local Langlands correspondence. Then
\begin{equation*}
   \gamma(s,\pi,r,\psi) = \left\{ \begin{array}{ll}
   						\gamma(s,{\rm Sym}^2\sigma,\psi) & \text{ if } {\bf G} = {\rm SO}_{2n+1} \\
						\gamma(s,\wedge^2 \sigma,\psi) & \text{ if } {\bf G} = {\rm SO}_{2n} \\
						\gamma(s,{}^\otimes{\rm I}(\sigma),\psi) & \text{ if } {\bf G} = {\rm U}_{2n}
				\end{array} . \right.
\end{equation*}
The $\gamma$-factors on the right hand side are those defined by Deligne and Langlands for Weil-Deligne representations.
}

\medskip

This is proved in \cite{hl11,hlarxiv} via a characterization of $\gamma$-factors in terms of local properties and their appearance in the global functional equation of partial $L$-functions. This is done with the help of a local-to-global result. The necessary properties for the uniqueness argument to hold are given in \S~6.5 below. The proof when ${\bf G} = {\rm GL}_{n_1+n_2}$ is given in the appendix. It is also possible to adapt the argument and include ${\bf G} = {\rm Sp}_{2n}$ and ${\rm U}_{2n+1}$, which are no longer part of the main induction since $r = r_1 \oplus r_2$.

\medskip

\noindent{\bf Theorem'.} \emph{Let $(E/F,\pi,\psi) \in \mathscr{L}({\bf M},{\bf G})$ and let $\sigma$ be as in the previous Theorem. If ${\bf G} = {\rm Sp}_{2n}$, $n>1$, then
\begin{equation*}
   C'_\psi(s,\pi_t,w_0) = \gamma(s,\sigma,\psi) \, \gamma(2s,\wedge^2\sigma,\psi).
\end{equation*}
If ${\bf G} = {\rm U}_{2n+1}$ and $\pi' = \pi \otimes \nu$, where $\nu$ is a character of $E^1$, extend $\nu$ to a smooth representation of ${\rm GL}_1(E)$ via Hilbert's theorem~90. Then
\begin{equation*}
   C'_\psi(s,\pi'_t,w_0) = \gamma_E(s,\pi \times \nu,\psi \circ {\rm Tr}_{E/F}) \, \gamma_F(2s,\pi \otimes \eta_{E/F},r_\mathcal{A},\psi).
\end{equation*}
}

\subsection{Properties of $\gamma$-factors} Local properties, ${\bf (i)} - {\bf (v)}$, needed to uniquely characterize generic $\gamma$-factors are now provided (see \cite{hl11,hlarxiv} and the appendix). The link to the global theory is done via the functional equation~${\bf (vi)}$.

\medskip

\noindent {\bf (i)} (Naturality). Let $(E/F,\pi,\psi) \in \mathscr{L}({\bf M},{\bf G})$. Let $\eta: E' \rightarrow E$ be an isomorphism such that $\eta_{F'}$ maps $F'$ into $F$. Then $\psi' = \psi \circ \eta \vert_{F'}$ is a non-trivial character of $F'$. Also, via $\eta$, $\pi$ defines a smooth irreducible representation $\pi'$ of ${\bf M}(F')$. Then $a_t = a_{t'}$ and
\begin{equation*}
   C'_\psi(s,\pi_t,w_0) = C'_\psi(s,\pi'_{t'},w_0).
\end{equation*}

\noindent {\bf (ii)} (Isomorphism). Let $\pi_1$ and $\pi_2$ be two irreducible representations $\pi_1$ and $\pi_2$ of $M$, generic with respect to the same non-degenerate character $\chi_M$, and such that $\pi_1 \simeq \pi_2$. Then $C_\chi(s,\pi_1,w_0) = C_\chi(s,\pi_2,w_0)$, see Theorem~3.1 of \cite{sha81}. This gives an isomorphism property
\begin{equation*}
   C'_\psi(s,(\pi_1)_{t_1},w_0) = C'_\psi(s,(\pi_2)_{t_2},w_0),
\end{equation*}
for $(E/F,\pi_i,\psi) \in \mathscr{L}({\bf M},{\bf G})$, $i = 1$, $2$, with $\pi_1 \simeq \pi_2$.

\medskip

\noindent {\bf Remark.} Because of this isomorphism property, it is possible to assume that, given $(E/F,\pi,\psi) \in \mathscr{L}({\bf M},{\bf G})$, the representation $\pi$ is $\psi_M$-generic.

\medskip

\noindent {\bf (iii)} Proposition~3.2 shows that local coefficients are compatible with class field theory.

\medskip

\noindent {\bf (iv)} Equation~\eqref{psidependence} explains the behavior of the local coefficient when the non-trivial character varies. 

\medskip

\noindent {\bf (v)} (Multiplicativity). First, assume that ${\bf G}$ is a classical group. Let $(E/F,\pi,\psi) \in \mathscr{L}({\bf M},{\bf G})$ be of degree $n$. Let $n = n_1 + \cdots + n_d$ and suppose $\pi$ is the generic constituent of
\begin{equation*}
   {\rm ind}_P^{{\rm GL}_n(E)} (\pi_1 \otimes \cdots \otimes \pi_d),
\end{equation*}
where ${\bf P} = {\bf M}{\bf N}$ is a parabolic subgroup of ${\rm GL}_n$ with Levi ${\bf M} = \prod_{k=1}^d {\rm GL}_{n_k}$. Then, multiplicativity gives a relationship for $\gamma$-factors involving of the inducing data
\begin{equation} \label{multgeneral}
   \gamma(s,\pi,r_i,\psi) = \prod_{j} \gamma(s,\pi_{i,j},r_{i,j},\psi),
\end{equation}
where the product ranges over a specific set of indices. For each $j$, $(E/F,\pi_{i,j},\psi) \in \mathscr{L}({\bf M}_j,{\bf G}_j)$, with ${\bf G}_j$ either a group of the same type as $\bf G$ or a general linear group.

Multiplicativity is stated explicitly in \cite{lomeli2009} when ${\bf G} = {\rm SO}_{2n+1}$, ${\rm Sp}_{2n}$, ${\rm SO}_{2n}$, and in \cite{hlarxiv} when ${\bf G} = {\rm U}_{2n}$. An explicit relationship for each of the remaining cases of ${\bf G} = {\rm GL}_{n_1+n_2}$ and ${\rm U}_{2n+1}$ is now given.

Now, assume that ${\bf G}$ is a general linear group. Let $(F,\pi,\psi) \in \mathscr{L}$ be of degree $(n_1,n_2)$, with $\pi = \pi_1 \otimes \tilde{\pi}_2$. Assume that $\pi_1$ is the generic constituent of
\begin{equation*}
   {\rm ind}_{P_1}^{{\rm GL}_{n_1}(F)}(\xi_1 \otimes \cdots \otimes \xi_d),
\end{equation*}
with each $\xi_i$ a representation of ${\rm GL}_{h_i}(F)$, $n_1 = h_1 + \cdots + h_d$. Also, assume that $\pi_2$ is the generic constituent of
\begin{equation*}
   {\rm ind}_{P_2}^{{\rm GL}_{n_2}(F)}(\tau_1 \otimes \cdots \otimes \tau_e),
\end{equation*}
with each $\tau_j$ a representation of ${\rm GL}_{m_j}(F)$, $n_2 = m_1 + \cdots + m_e$. Then
\begin{equation} \label{multGLn}
   \gamma(s,\pi,\psi) = \prod_{k,l} \gamma(s,\xi_k \otimes \tilde{\tau}_l,\psi),
\end{equation}
where the product ranges over all $k$, $1 \leq k \leq d$, and all $l$, $1 \leq l \leq e$. Equation~\eqref{multgeneral} can now be explicitly interpreted in this case with $\pi_{i,j}$ of the form $\xi_k \otimes \tilde{\tau}_l$ and $r_{i,j}$ the appropriate standard representation.

Next, let $(E/F,\pi,\psi) \in \mathscr{L}({\bf M},{\rm U}_{2n+1})$, and suppose $\pi = \pi' \otimes \nu'$ is $\psi_M$-generic. Extend $\nu'$ to a character $\nu$ of ${\rm GL}_1(E)$. Let $n = n_1 + \cdots + n_d$ and suppose $\pi'$ is the $\psi$-generic constituent of
\begin{equation*}
   {\rm ind}_{P}^{{\rm GL}_n(E)}(\pi_1 \otimes \cdots \otimes \pi_d),
\end{equation*}
where ${\bf P} = {\bf M} {\bf N}$ is a parabolic subgroup of ${\rm GL}_n$ with Levi ${\bf M} = \prod_{i=1}^{d} {\rm GL}_{n_i}$ and, for each $i$, $\pi_i$ is a supercuspidal $\psi$-generic representation of ${\bf GL}_{n_i}(E)$. Then
\begin{align*}
   C'_\psi(s,\pi' \otimes \nu',w_0) 	& = \prod_{i=1}^{d} C'_\psi(s,\pi_i \otimes \nu,w_0^i)  \\
   				& \times \prod_{i<j} C_{\psi}(2s,\pi_i \otimes \tilde{\pi}_j^{\rm conj.},w_0^{ij}).
\end{align*}
Here, $C'_{\psi}(s,\pi_i \otimes \nu,w_0^i)$ is a local coefficient of ${\rm U}_{2n_i+1}$ and $C_\psi(2s,\pi_i \otimes \tilde{\pi}_j^{\rm conj.},w_0^{ij})$ is a local coefficient of ${\rm GL}_{n_i + n_j}(E)$ equal to $\gamma(2s,\pi_i \times \pi_j^{\rm conj.},\psi_E)$. By grouping the factors involving $s$ and $2s$ together, equation~\eqref{multgeneral} takes the following explicit form for each of the two corresponding $\gamma$-factors:
\begin{subequations}
\begin{align}
   \gamma(s,\pi,r_1,\psi)	& = \gamma_E(s,\pi' \times \nu,\psi_E) = \prod_{k=1}^{d} \gamma_E(s,\pi_k \times \nu,\psi_E) \\
   \gamma(s,\pi,r_2,\psi)	& = \prod_{k=1}^d \gamma(s,\pi_k \otimes \eta_{E/F},r_\mathcal{A},\psi) \prod_{k<l} \gamma_E(s,\pi_k \times \pi_l^{\rm conj},\psi_E).
\end{align}
\end{subequations}

\medskip

\noindent{\bf (vi)} (Global functional equation). With the notation of \S~5, ${\bf G}$ is defined over a global function field $k$. It is possible to embed ${\bf G}$ into a group $\widetilde{\bf G}$ sharing the same derived group as ${\bf G}$, and satisfying $H^1(Z_{\widetilde{G}}) = \left\{ 1 \right\}$. Let $\widetilde{\bf B} = \widetilde{\bf T}{\bf U}$ be a Borel subgroup of $\widetilde{\bf G}$ defined over $k$. Then $\widetilde{\bf T}(k)$ acts transitively on generic characters of ${\bf U}(k) \backslash {\bf U}(\mathbb{A}_k)$. Hence, it is possible to find a $t \in \widetilde{\bf T}(k)$ such that $\pi_t$ is globally $\psi$-generic, where $\pi_t(m) = \pi(t^{-1}mt)$. See the appendix to \cite{ckpss}, where the homological discussion holds for a global function field. The functional equation of \S~6, now becomes one involving $\gamma$-factors
\begin{equation*}
   L^S(s,\pi,r) = \prod_{v \in S} \gamma(s,\pi_v,r_v) L^S(1-s,\tilde{\pi},r),
\end{equation*}
where $r = \rho_{n_1} \otimes \tilde{\rho}_{n_2}$, ${\rm Sym}^2 \rho_n$, $\wedge^2 \rho_n$, or $r_\mathcal{A}$. If ${\bf G} = {\rm Sp}_{2n}$ (with $n>1$) or ${\rm U}_{2n+1}$, there are two individual functional equations:
\begin{equation*}
   L^S(s,\pi,r_i) = \prod_{v \in S} \gamma(s,\pi_v,r_{i,v}) L^S(1-s,\tilde{\pi},r_i),
\end{equation*}
for $i = 1$, $2$, since $r = r_1 \oplus r_2$ in the notation of \S~5.2. Notice that, if ${\bf G} = {\rm Sp}_2 = {\rm SL}_2$, one retrieves the theory of $L$-functions for an automorphic representation of ${\rm GL}_1$, i.e., for $\pi = \chi$ a Gr\"ossencharakter of ${\rm GL}_1(\mathbb{A}_k)$.

\subsection{Additional properties of $\gamma$-factors} \texttt{ }

\medskip

\noindent{\bf (vii)} (Twists by unramified characters) Let $(E/F,\pi,\psi) \in \mathscr{L}({\bf M},{\bf G})$. If ${\bf G} = {\rm GL}_{n_1+n_2}$ and $\pi = \pi_1 \otimes \tilde{\pi}_2$, then
\begin{equation*}
   \gamma(s + s_0,\pi_1 \otimes \tilde{\pi}_2,\psi) = \gamma(s,(\left| \det(\cdot) \right|_F^{s_0} \pi_1) \otimes \tilde{\pi}_2,\psi).
\end{equation*}
If ${\bf G} = {\rm SO}_{2n}$, ${\rm SO}_{2n+1}$ or ${\rm U}_{2n}$, then
\begin{equation*}
   \gamma(s + s_0,\pi,r,\psi) = \gamma(s,\left| \det(\cdot) \right|_F^{\frac{s_0}{2}} \pi, r,\psi).
\end{equation*}

Notice that, when ${\rm Sp}_{2n}$ and ${\rm U}_{2n+1}$, the relationship corresponding to twists by unramified characters can be obtained from these cases for each of the $\gamma$-factors $\gamma(s,\pi,r_i,\psi)$, $i=1$ or $2$.

\medskip

\noindent{\bf (viii)} (Local functional equation) Let $(E/F,\pi,\psi) \in \mathscr{L}({\bf M},{\bf G})$, then
\begin{equation*}
   \gamma(s,\pi,r_i,\psi) \, \gamma(1-s,\tilde{\pi},r_i,\overline{\psi}) = 1.
\end{equation*}

\section{$L$-functions, $\varepsilon$-factors, and their global functional equation.}

Local $L$-functions and root numbers are first defined for tempered representations. Then in general with the aid of Langlands classification for $\mathfrak{p}$-adic groups combined with multiplicativity and twists by unramified characters. When $\bf G$ is a general linear group, the approach taken here is in accordance with the theory for principal $L$-functions \cite{j} and Rankin-Selberg convolutions \cite{jpss}.

\subsection{$L$-functions and tempered representations} With the notation of \S~6.1, given $(E/F,\pi,\psi) \in \mathscr{L}({\bf M},{\bf G})$, call it tempered (resp. supercuspidal) if the irreducible representation $\pi$ of $M$ is tempered (resp. supercuspidal). For tempered $(E/F,\pi,\psi)$, let $P_{\pi,i}(t)$ be the unique polynomial satisfying $P_{\pi,i}(0) = 1$ such that $P_{\pi,i}(q^{-s})$ is the numerator of $\gamma(s,\pi,r_i,\psi)$, $1 \leq i \leq m_r$. Define
\begin{equation}
   L(s,\pi,r_i) = P_{\pi,i}(q^{-s})^{-1}.
\end{equation}

Notice that the zeros of a local coefficient are poles of the corresponding intertwining operator ${\rm A}(s,\pi,w_0)$. Hence, the definition of a local $L$-function does not depend on the character $\psi$.

\medskip

\subsection{Proposition} \emph{Let $(E/F,\pi,\psi) \in \mathscr{L}({\bf M},{\bf G})$ be tempered. Then $L(s,\pi,r_i)$, $1 \leq i \leq m_r$, is holomorphic for $\Re(s)>0$.}

\medskip

\noindent {\bf Proof.} First, consider the cases ${\bf G} = {\rm GL}_{n_1+n_2}$, ${\rm SO}_{2n+1}$, ${\rm SO}_{2n}$ or ${\rm U}_{2n}$. From the definitions, the zeros of $C_\psi(s,\pi,w_0)$ are those of $P_{\pi,1}(q^{-s})$. The proposition follows in these cases, since it is known that the intertwining operators are holomorphic for $\Re(s)>0$ when $\pi$ is tempered (see Theorem~5.3.5.4 of \cite{sil}; also, see Proposition~IV.2.1 of \cite{w}). Next, consider the cases ${\bf G} = {\rm Sp}_{2n}$ or ${\rm U}_{2n+1}$. Then, $L(s,\pi,r_1)$ arises from a local coefficient of ${\rm GL}_{n+1}$, and $L(2s,\pi,r_2)$ arises arises from a local coefficient of either ${\rm SO}_{2n}$ or ${\rm U}_{2n}$, each of which is holomorphic for $\Re(s)>0$ by the previous cases.

\medskip

\subsection{Tempered representations of ${\rm GL}_n$} Every tempered representation of ${\rm GL}_n(F)$ arises as the unique irreducible submodule of a parabolically induced representation
\begin{equation*}
   {\rm ind}_P^{{\rm GL}_n(F)} (\delta_1 \otimes \cdots \otimes \delta_r),
\end{equation*}
where each $\delta_i$ is a discrete series representation. This is part of Theorem~9.7 of \cite{z}, which is stated for quasi-tempered representations, i.e., tempered representations twisted by a suitable character.

In order to be more precise, given a supercuspidal representation $\rho$ of ${\rm GL}_d(F)$ and a positive integer $a$, let $\delta(\rho,a)$ be the irreducible submodule of
\begin{equation*}
   {\rm ind}_{P_a}^{{\rm GL}_d(F)}(\nu^{\frac{a-1}{2}} \rho \otimes \cdots \otimes \nu^{-\frac{a-1}{2}} \rho),
\end{equation*}
where $\nu = \left| \det(\cdot) \right|_F$ and ${\bf P}_a = {\bf M}_a {\bf N}_a$ is the parabolic with Levi ${\bf M}_a = \prod_{i=1}^{a} {\rm GL}_d$. If $\rho$ is unitary, then $\delta(\rho,a)$ is a unitary discrete series representation of ${\rm GL}_{ad}(F)$. All discrete series representations of ${\rm GL}_n(F)$ are of this form.

Notice that, if $\pi_i = \delta(\rho_i,a_i)$, $i =1$ or $2$, are discrete series representations of ${\rm GL}_{n_i}(F)$, then $L(s,\pi_1 \times \pi_2)$ is holomorphic for $\Re(s) > 0$. Similarly for $L(s,\pi,r)$, with $\pi = \delta(\rho,a)$ in the previous proposition.

\subsection{Root numbers and tempered representations.} Let $(E/F,\pi,\psi) \in \mathscr{L}({\bf M},{\bf G})$ be tempered. Then, there is a monomial $\varepsilon(s,\pi,r_i,\psi)$ in the variable $Z = q^{-s}$, called the local root number or $\varepsilon$-factor, given by
\begin{equation}
   \varepsilon(s,\pi,r_i,\psi) = \gamma(s,\pi,r_i,\psi) \dfrac{L(s,\pi,r_i)}{L(1-s,\tilde{\pi},r_i)}.
\end{equation}

\medskip

That the root number is indeed a monomial in $Z = q^{-s}$, follows from Proposition~7.2 and the local functional equation of $\gamma$-factors~6.6(viii).

The multiplicativity property of $\gamma$-factors is inherited for local $L$-functions and $\varepsilon$-factors when $(E/F,\pi,\psi) \in \mathscr{L}({\bf M},{\bf G})$ is tempered. With the notation of equation~\eqref{multgeneral}, multiplicativity of local factors reads
\begin{subequations}
\begin{align}
   L(s,\pi,r_i) & = \prod_j L(s,\pi_{i,j},r_{i,j}) \\
   \varepsilon(s,\pi,r_i,\psi) & = \prod_{j} \varepsilon(s,\pi_{i,j},r_{i,j},\psi).
\end{align}
\end{subequations}

\subsection{Local factors in general} It is now possible to proceed to the general case. For each $i$, $1 \leq i \leq d$, let $\tau_{i,0}$ be a tempered representation of ${\rm GL}_{n_i}(F)$. Take a sequence of real numbers $t_1 > \cdots > t_d$ and let $\tau_i = \left| \det(\cdot) \right|_F^{t_i} \tau_{i,0}$, for each $i$. Let $n = n_1 + \cdots + n_d$ and let ${\bf P}' = {\bf M}'{\bf N}'$ be the parabolic subgroup of ${\rm GL}_n$, with Levi ${\bf M}' = \prod_{i=1}^{d} {\rm GL}_{n_i}$. Then, the unitarily induced representation
\begin{equation} \label{generalinduced}
   \tau = {\rm ind}_{P'}^{{\rm GL}_n(F)}(\tau_1 \otimes \cdots \otimes \tau_d),
\end{equation}
has a unique Langlands' quotient representation $\pi = {\rm J}(\tau)$. Every irreducible representation $\pi$ of ${\rm GL}_n$ is of this form.

Now, with the notation of \S~6.1 when $\bf G$ is a general linear group, let $(F,\pi,\psi) \in \mathscr{L}$ be of degree $(n_1,n_2)$, with $\pi = \pi_1 \otimes \tilde{\pi}_2$. It is then possible to assume that $\pi_1 = {\rm J}(\xi)$ and $\pi_2 = {\rm J}(\tau)$, where $\xi = {\rm ind}_{P_{n_1}}^{{\rm GL}_{n_1}(F)}(\xi_1 \otimes \cdots \otimes \xi_d)$ and $\rho = {\rm ind}_{P_{n_2}}^{{\rm GL}_{n_2}(F)}(\tau_1 \otimes \cdots \otimes \tau_e)$ are as in equation~\eqref{generalinduced} with $\xi_k = \left| \det(\cdot) \right|_F^{u_k} \xi_{k,0}$ and $\tau_l = \left| \det(\cdot) \right|_F^{v_l} \tau_{l,0}$. Equation~\eqref{multGLn} combined with twists by unramified characters gives
\begin{equation} \label{multtempGLn}
   \gamma(s,\pi_1 \otimes \tilde{\pi}_2,\psi) = \prod_{k,l} \gamma(s + u_k + v_l ,\xi_{k,0} \otimes \tilde{\tau}_{l,0},\psi).
\end{equation}
Here, each $\gamma$-factor in the product corresponds to tempered representations. This allows one to define local $L$-functions and root numbers in general by setting
\begin{subequations}
\begin{align} \label{LrootdefGLn}
   L(s,\pi) & = \prod_{k,l} L(s+u_k+v_l,\xi_{k,0} \otimes \tilde{\tau}_{l,0}) \\
   \varepsilon(s,\pi,\psi) & = \prod_{k,l} \varepsilon(s+u_k+v_l,\xi_{k,0} \otimes \tilde{\tau}_{l,0},\psi).
\end{align}
\end{subequations}

Next, when ${\bf G}$ is a classical group, let $(E/F,\pi,\psi) \in \mathscr{L}({\bf M},{\bf G})$ be such that $\pi = {\rm J}(\tau)$ as in equation~\eqref{generalinduced}. Then local $L$-functions and root numbers are defined using multiplicativity, equation~\eqref{multgeneral}, and twists by unramified characters, property~6.6(vii). It is now possible to obtain explicit formulas for each of the classical groups in terms of tempered inducing data. In particular, Asai local factors are defined in general in \S~4.3 of \cite{hlarxiv}. It is now an exercise to explicitly do this for exterior square and symmetric square local factors.

\subsection{Global functional equation} With the notation of \S~5, let $\pi$ be a cuspidal automorphic representation of ${\bf M}(\mathbb{A}_k)$. Local $L$-functions and root numbers are now defined at every place of $k$. Thus, it is possible to define the completed $L$-functions:
\begin{equation*}
   L(s,\pi,r_i) = \prod_v L(s,\pi_v,r_{i,v}) \ \text{ and } \ \varepsilon(s,\pi,r_i) = \prod_v \varepsilon(s,\pi_v,r_{i,v},\psi_v),
\end{equation*}
for each $i$, $1 \leq i \leq m_r$. The final form of the functional equation is now obtained from Property~6.5(vi) for each $i$:
\begin{equation*}
   L(s,\pi,r_i) = \varepsilon(s,\pi,r_i) L(1-s,\tilde{\pi},r_i).
\end{equation*}
Notice that $\varepsilon(s,\pi_v,r_{i,v},\psi_v)$ is trivial outside a finite set $S$ of places of $k$ and that the global root number does not depend on the character $\psi$. Furthermore, rationality of completed $L$-functions follows from Corollary~5.4. To summarize, we have completed a linearly ordered proof of the following theorem.

\subsection{Theorem} \emph{Let $k$ be a global function field with finite field of constants $\mathbb{F}_q$. If ${\bf G}$ is a split classical group or a quasi-split unitary group, let ${\bf M}$ be a Siegel Levi subgroup; if ${\bf G}$ is a general linear group let ${\bf M}$ be a maximal Levi subgroup. Let $\pi$ be a cuspidal automorphic representation of ${\bf M}(\mathbb{A}_k)$. Then, the $L$-functions $L(s,\pi,r_i)$ are ${\rm nice}$, i.e., they satisfy:}
\begin{itemize}
   \item[(i)] (Rationality) \emph{$L(s,\pi,r_i)$ has a meromorphic continuation to a rational functions on $q^{-s}$;}
   \item[(ii)] (Functional equation) $L(s,\pi,r_i) = \varepsilon(s,\pi,r_i) L(1-s,\tilde{\pi},r_i)$.
\end{itemize}

\bigskip

{\sc \Small Luis Lomel\'i}

\emph{\Small E-mail address: }\texttt{\Small lomeli@math.purdue.edu}

\emph{\Small E-mail address: }\texttt{\Small lomeli@math.ou.edu}

\end{document}